\documentclass[11pt]{article}

\usepackage{a4wide}
\usepackage{amsfonts} 
\usepackage{amsmath}
\usepackage{amssymb}
\usepackage{mathrsfs}
\usepackage{pstricks}
\usepackage{pst-plot}
\usepackage[all]{xy}
\usepackage{float}

\newtheorem{Def}{Definition}[section]
\newtheorem{Lem}[Def]{Lemma}
\newtheorem{Thm}[Def]{Theorem}
\newtheorem{Prop}[Def]{Proposition}
\newtheorem{Cor}[Def]{Corollary}
\newtheorem{Remark}[Def]{Remark}
\newenvironment{Rem}{\begin{Remark}\rm}{\end{Remark}}

\newcommand{\cst}{\mathrm{C}^*\!}
\newcommand{\RR}{{\mathbb R}}
\newcommand{\CC}{{\mathbb C}}
\newcommand{\NN}{{\mathbb N}}
\newcommand{\ZZ}{{\mathbb Z}}
\newcommand{\TT}{{\mathbb T}}

\newcommand{\id}{\mathrm{id}}
\newcommand{\comp}{\!\circ\!}
\newcommand{\bez}{\setminus}
\newcommand{\la}{\langle}
\newcommand{\ra}{\rangle}

\newcommand{\ort}{\perp}
\newcommand{\Tr}[1]{\mathrm{Tr}\!\left(#1\right)}

\newcommand{\its}[3]{\left(#1\,\vline\,#2\,\vline\,#3\right)}
\newcommand{\K}[1]{{\mathcal K}\!\left(#1\right)}
\newcommand{\B}[1]{B\!\left(#1\right)}
\newcommand{\eps}{\varepsilon}
\newcommand{\ph}{\varphi}
\newcommand{\Sch}{{\mathscr S}}
\newcommand{\del}{\delta}
\newcommand{\UA}{\mathfrak{A}}
\newcommand{\UB}{\mathfrak{B}}
\newcommand{\UC}{\mathfrak{C}}
\newcommand{\Prod}{\prod\limits}
\newcommand{\Sum}{\sum\limits}

\newcommand{\SUM}[2]{\displaystyle\sum\limits_{#1}^{#2}}

\newcommand{\Int}{\int\limits}
\newcommand{\INt}[1]{\displaystyle\int\limits_{#1}}
\newcommand{\INT}[2]{\displaystyle\int\limits_{#1}^{#2}}
\newcommand{\ulamek}[2]{{\textstyle\frac{#1}{#2}}}
\newcommand{\tens}{\otimes}
\newcommand{\Aut}[1]{\mathrm{Aut}\!\left(#1\right)}
\newcommand{\Mor}[2]{\mathrm{Mor}\!\left(#1,#2\right)}
\newcommand{\M}[1]{M\!\left(#1\right)}
\newcommand{\aff}{\,\,\!\eta\,\,\!}

\newcommand{\Rep}[2]{\mathrm{Rep}\!\left(#1,#2\right)}
\newcommand{\dplus}{\,\dot{+}\,}
\newcommand{\D}[1]{D\!\left(#1\right)}
\newcommand{\spec}[1]{\mathrm{Sp}\,#1}
\newcommand{\Spec}[1]{\mathrm{Sp}\!\left(#1\right)}
\newcommand{\phase}[1]{\mathrm{Phase}\,#1}
\newcommand{\Phase}[1]{\left(\mathrm{Phase}\,#1\right)}
\newcommand{\PHase}[1]{\mathrm{Phase}\left(#1\right)}
\newcommand{\PHAse}[1]{\,\mathrm{Phase}\left(#1\right)}

\newcommand{\Zrel}[3]{\xymatrix{{#1}\ar@{-o}[r]^-{\scriptscriptstyle #3}&{#2}}}

\newcommand{\Exp}[1]{\exp{\left(#1\right)}}
\newcommand{\re}[1]{\mathrm{Re}\,#1}
\renewcommand{\Re}[1]{\mathrm{Re}\left(#1\right)}
\newcommand{\Rre}[1]{\,\Re{#1}}
\newcommand{\rre}[1]{\,\re{#1}}
\newcommand{\im}[1]{\mathrm{Im}\,#1}
\renewcommand{\Im}[1]{\mathrm{Im}\left(#1\right)}
\newcommand{\Iim}[1]{\,\Im{#1}}
\newcommand{\iim}[1]{\,\im{#1}}
\newcommand{\Ci}[1]{C_{\infty}\!\left(#1\right)}
\newcommand{\Cb}[1]{C_{\mathrm{bounded}}\!\left(#1\right)}

\newcommand{\refeq}[1]{{\rm (\ref{#1})}}
\newcommand{\proof}{\noindent{\sc Proof.~}}
\newcommand{\qed}{\hfill{\sc Q.E.D.}\medskip}
\newcommand{\An}[2]{#1\!\left[#2\right]}
\newcommand{\an}[3]{#1\!\left[#2\right]\left(#3\right)\!}
\newcommand{\rh}{\boldsymbol{\rho}\:\!}
\newcommand{\rhinv}{\rh^{-1}}
\newcommand{\rhbar}{\overline{\rh}}
\newcommand{\Yq}{{\Gamma_{\!\scriptscriptstyle q}}}
\newcommand{\Ybar}{\overline{\Gamma}_{\!\scriptscriptstyle q}}

\newcommand{\Cq}{C_{\!\scriptscriptstyle q}}
\newcommand{\y}{\gamma}
\newcommand{\qbar}{\overline{q}}
\newcommand{\Fq}{{\mathbb F}_{\!\scriptscriptstyle q}\:\!}
\newcommand{\Hq}{{\mathbb H}_{\scriptscriptstyle q}\:\!}
\newcommand{\Phiq}{{\Phi}_{\!\scriptscriptstyle q}\:\!}
\newcommand{\Dd}{{\mathscr D}}
\newcommand{\Gg}{{\mathscr G}}
\newcommand{\PH}{\mathrm{Ph}}

\newcommand{\pH}[1]{\PH\!\left(#1\right)}
\newcommand{\Mu}{\mu}
\newcommand{\ahat}{\widehat{a}}
\newcommand{\bhat}{\widehat{b}}
\newcommand{\Vtil}{\widetilde{V}}
\newcommand{\Wtil}{\widetilde{W}}
\newcommand{\Qhat}{\widehat{Q}}
\newcommand{\Hh}{{\mathscr H}}
\newcommand{\xbar}{\overline{x}}
\newcommand{\zbar}{\overline{z}}
\newcommand{\void}[1]{}

\newcommand{\Section}[1]{\setcounter{equation}{0}\section{#1}}

\begin{document}

\title{New quantum ``$az+b$'' groups}
\author{Piotr Miko\l{}aj So\l{}tan\thanks{
Research partially supported by Komitet Bada\'n Naukowych grant 
no.~2PO3A04022, the Foundation for Polish Science and Deutsche 
Forschungsgemeinschaft.} \\
\small Department of Mathematical Methods in Physics,\\
\small Faculty of Physics, University of Warsaw\\
\small\tt piotr.soltan@fuw.edu.pl}
\maketitle

\abstract{We construct quantum ``$az+b$'' groups for new values of the
deformation parameter. Along the way we introduce new special functions and
study their analytic properties as well as analyze the commutation relations
determined by the choice of parameter.}


{\tiny
\tableofcontents}

\section{Introduction}

The aim of this paper is to carry out the construction of quantum ``$az+b$''
presented in \cite{azb}, for new values of the deformation parameter. The
construction procedure form \cite{azb} will be repeated with necessary
modifications. The main problems lie in developing the machinery of special
functions needed for the construction and applying this machinery to analysis
of appropriate commutation relations. We also propose ways to streamline
some aspects of the construction known from \cite{azb}. 

\subsection{Algebraic quantum ``$az+b$'' groups}\label{alg_azb}

On the level of Hopf $*$-algebra the quantum ``$az+b$'' groups we will be
interested in are described as follows: let $\lambda$ be a non zero complex
number. Let $\Hh$ be a unital $*$-algebra generated by three normal elements 
$a$, $a^{-1}$ and $b$ with the relations
\[
\begin{array}{r@{\;=\;}l@{\qquad}r@{\;=\;}l@{\smallskip}}
a^{-1}a&I,&aa^{-1}&I,\\
ab&\lambda ba,&ab^*&b^*a.
\end{array}
\]  
Then $\Hh$ can be given a structure of a Hopf $*$-algebra by
\[
\begin{array}{r@{\;=\;}l@{\smallskip}}
\del(a)&a\tens a,\\
\del(b)&a\tens b+b\tens a.
\end{array}
\]
The coinverse and counit are given on generators in the following way:
\[
\begin{array}{r@{\;=\;}l@{\qquad}r@{\;=\;}l@{\smallskip}}
\kappa(a)&a^{-1},&\epsilon(a)&1,\\
\kappa(b)&-a^{-1}b,&\epsilon(b)&0.
\end{array}
\]

An important fact about $\Hh$ is that the coinverse has a {\em polar
decomposition} (cf.~\cite[Proposition 2.4]{VDH}. By existence of a polar 
decomposition of $\kappa$ we mean that
there exist a $*$-antiautomorphism $R$ of $\Hh$ and a one parameter group
$(\tau_t)_{t\in\RR}$ of $*$-automorphisms of $\Hh$ such that for any linear
functional $f$ on $\Hh$ and any $x\in\Hh$ the map
\[
\RR\ni t\longmapsto f\bigl(\tau_t(x)\bigr)\in\CC
\]
has an extension to an entire function and $\kappa=R\comp\tau_{\frac{i}{2}}$,
where $\tau_{\frac{i}{2}}$ is an automorphism of $\Hh$ obtained by analytic
continuation of the group $(\tau_t)_{t\in\RR}$. 

For $\lambda=e^{\frac{2\pi i}{n}}$ we can impose further condition that $a^n$
and $b^n$ be self adjoint. This corresponds to taking a quotient of $\Hh$ by
the ideal generated by $a^n-(a^n)^*$ and $b^n-(b^n)^*$. It turns out that this
is a Hopf ideal, i.e.~the quotient inherits Hopf $*$-algebra structure.
Moreover the coinverse on the resulting Hopf $*$-algebra still has a polar
decomposition. 

\subsection{Quantum ``$az+b$'' groups of S.L.~Woronowicz}\label{SLWazb}

A great challenge taken up in \cite{azb} was to construct the quantum
``$az+b$'' group on $\cst$-algebra level. The construction was based on the
Hopf $*$-algebraic picture with $\lambda=e^{\frac{2\pi i}{n}}$ with $n\geq3$ and
the relations $a^n-(a^n)^*$ and $b^n-(b^n)^*$ were translated into {\em 
spectral conditions}\/ for $a$ and $b$. On $\cst$-algebra level the generators
$a$ and $b$ become unbounded operators. The condition that their
$n^{\text{\tiny th}}$ powers be self adjoint are equivalent to the condition
that their spectra be contained in the closure of the set
\[
\bigcup_{k=0}^{2n-1}e^{\frac{k\pi i}{n}}\RR_+.
\] 
It was noticed that this set is a multiplicative (and self dual) 
subgroup of $\CC\bez\{0\}$
and the commutation relations between $a$ and $b$ were written in the {\em
Weyl form}\/ with use of a bicharacter on this group (cf.~also Subsection
\ref{komrel}).

This very successful construction was then performed in a different setting.
In \cite[Appendix B]{azb} a construction of quantum ``$az+b$'' groups was also
done for $0<\lambda<1$. One could no longer take the Hopf $*$-algebra with
$a^n$ and $b^n$ self adjoint. Therefore one starts with the Hopf $*$-algebra
$\Hh$ (see the previous subsection) and then imposes the condition that the
spectra of $a$ and $b$ be contained in the closure of the following
multiplicative (and self dual) subgroup of $\CC\bez\{0\}$:
\[
\bigl\{z\in\CC:\:|z|\in\lambda^{\frac{1}{2}\ZZ}\bigr\}.
\] 
It turns out that the construction works and we obtain quantum ``$az+b$''
groups for real deformation parameter. Note that this spectral condition does
not correspond to any quotient of $\Hh$. 

The quantum ``$az+b$'' groups we aim to construct in this paper will be of
similar nature. They do not correspond to any quotient of the Hopf $*$-algebra
$\Hh$. Nevertheless some aspects of these groups are similar to those with
deformation parameter assuming the value of an even root of unity.

\subsection{Description of the paper}

Let us now briefly describe the contents of the paper. In Section
\ref{SpecFun} we describe the values of the deformation parameter which we
will use throughout the paper. Then we introduce the subsets
$\Yq,\Ybar\subset\CC$ and define three special functions on these sets with
values in the unit circle $\TT$. We provide detailed proofs of special
relations between these functions, discuss their asymptotic behaviour and
uniqueness properties.

In Section \ref{opeq} we define and analyze the commutation relations needed
for the construction of new quantum ``$az+b$'' groups. This section lays
foundations for the study of our quantum groups on $\cst$-algebra level.

Section \ref{affrel} is devoted to the algebraic consequences of commutation
relations studied in Section \ref{opeq}. Since we are dealing with unbounded
operators, all results are formulated with use of the affiliation relation for
$\cst$-algebras. Some results of \cite{azb} are generalized in such a way that
they can be applied to our construction. This is accomplished mainly through
unveiling the general mechanisms behind them. The concept of a $\cst$-algebra
generated by unbounded elements as well as by a quantum family of affiliated
elements plays an important role in these considerations. 

In Section \ref{MU} we define and study the basic object leading to the
construction of our quantum ``$az+b$'' group, i.e~the multiplicative unitary
operator. Using the machinery of operator equalities and special functions
developed in Sections \ref{opeq} and \ref{SpecFun} we show that this
multiplicative unitary is {\em modular}\/ in the sense of \cite{modmu}.

In the last section we identify the $\cst$-algebra of continuous functions
vanishing at infinity on our deformation of ``$az+b$'' group. Then we describe
the generators $a$ and $b$ and introduce the quantum group structure. It is
all done in accordance with the general procedures known e.g.~from \cite{mu}.
Then using the latest results of S.L.~Woronowicz (\cite{haar}) we introduce
the right Haar measure on our quantum group. Thus we show that the constructed
object falls into the category studied in \cite{mnw,kv}.

\Section{Three special functions and their properties}\label{SpecFun}

The construction of new quantum ``$az+b$'' groups on $\cst$-algebra level in
Section \ref{NewQG} will extend the framework presented in \cite{azb}. It will
be governed by a deformation parameter $q$. We shall impose certain conditions 
on the value of this parameter. The admissible values of the deformation
parameter are
\[
q=\Exp{{\rhinv}},
\]
where $\rh$ is a complex number such that $\re{\rh}<0$ and 
$\im{\rh}=\ulamek{N}{2\pi}$ with $N$ an even natural number. The number 
$\rh$ or the pair $(N,\re{\rh})$ can equally be taken to be the fundamental
parameters of our theory. The choice of $q$ as the main parameter is motivated
by the traditional formulation of commutation relations discussed in Subsections
\ref{alg_azb} and \ref{SLWazb} and Section \ref{opeq}. The choice of $\rh$
determines a choice of logarithm of $q$: for any $z\in\CC$ we set
\[
q^z=\Exp{\ulamek{z}{\rh}}.
\]

All our constructions work equally well for negative values of $N$. The
restriction to the case $N\in 2\NN$ helps keep our notation simpler. The case
of all admissible values of $q$ (including the ones with $N<0$) is treated in
detail in \cite{phd}.

The special functions discussed in this section will be defined on subsets of
$\CC$ defined with use f the parameter $q$. Let $\Yq$ be the multiplicative
subgroup of $\CC\bez\{0\}$ generated by $q$ and $\{q^{it}:\:t\in\RR\}$ and let
$\Ybar$ be the closure of $\Yq$ in $\CC$.
 
\subsection{The bicharacter $\chi$}

The group $\Yq$ with the topology inherited from $\CC\bez\{0\}$ is isomorphic
as a locally compact group to $\ZZ_N\times\RR$. It is therefore self dual. The
isomorphism of $\Yq$ to its dual group can be encoded by a non degenerate
bicharacter defined on $\Yq\times\Yq$. We shall choose a particular
bicharacter and use it throughout the paper.

\begin{Prop}
There exists a unique continuous function $\chi\colon\Yq\times\Yq\to\TT$ such
that for all $\y,\y'\,y''\in\Yq$ 
\[
\begin{array}{r@{\;=\;}l@{\smallskip}}
\chi(\y,\y')&\chi(\y',\y),\\
\chi(\y,\y')\chi(\y,\y'')&\chi(\y,\y'\y'')
\end{array}
\]
and
\begin{equation}\label{bchar}
\begin{array}{r@{\;=\;}l@{\smallskip}}
\chi(\y,q)&\phase{\y},\\
\chi(\y,q^{it})&|\y|^{it}
\end{array}
\end{equation}
for all $\y\in\Yq$, $t\in\RR$. The function $\chi$ is non degenerate,
i.e.~$\chi(\y,\y')=1$ for all $\y'\in\Yq$ implies $\y=1$.
\end{Prop} 

For $\y=q^{n}q^{it}$ and $\y'=q^{n'}q^{it'}$ we have
\begin{equation}\label{DefChi}
\chi(\y,\y')=e^{i(nn'-tt')\im{\rhinv}}e^{i(nt'+n't)\re{\rhinv}}.
\end{equation}
The expression of $\y$ as $q^{n}q^{it}$ is not unique. In fact we have
\[
q^{n}q^{it}=q^{n+N}q^{i\left(t-\frac{\re{\rh}}{\im{\rh}}\right)}
\]
for all $n\in\ZZ$ and $t\in\RR$. Nevertheless one can easily check that the
definition \refeq{DefChi} is does not depend on the way $\y$ and $\y'$ are 
expressed in the form $q^{n}q^{it}$ and $q^{n'}q^{it'}$ respectively.  

\subsection{The Fresnel function $\alpha$}

\begin{Prop}\label{alfa}
There exists a continuous function $\alpha\colon\Yq\to\TT$ such that for all
$\y,\y'\in\Yq$
\begin{equation}\label{AlfaChi}
\chi(\y,\y')=\frac{\alpha(\y\y')}{\alpha(\y)\alpha(\y')}.
\end{equation}
and 
\[
\alpha(\y)=\alpha(\y^{-1}).
\]
\end{Prop}

\proof
The function $\alpha$ is unique up to multiplication by a $\ZZ_2$-valued
character of $\Yq$. The formula
\begin{equation}\label{Alfa}
\alpha(q^{n}q^{it})=\phase{\frac{(n+it)^2}{2\rh}}
=\Exp{i\iim{\frac{(n+it)^2}{2\rh}}}
\end{equation}
defines a function with required properties.
\qed

We chose a particular formula for $\alpha$ only to make the exposition more
transparent. The important properties are those described in the statement of
Proposition \ref{alfa}. Formula \refeq{Alfa} will make it easier to proceed
with computations, but whenever we use it only the absolute value of $\alpha$
(which is uniquely determined) will enter our considerations. 

The pair of functions $(\chi,\alpha)$ is in many aspects analogous to the pair
of functions $\RR\times\RR\ni(p,x)\mapsto e^{ipx}\in\TT$ and
$\RR\ni x\mapsto e^{i\frac{x^2}{2}}\in\TT$. The latter one enters the formula 
for the Fresnel integral from which we borrow the name {\em Fresnel function}\/ 
for $\alpha$.

\subsection{The quantum exponential function $\Fq$}\label{qefun}

From the fact that $N$ is even it follows that $-1\in\Yq$ and consequently
$-q^{-2k}$ belongs to $\Yq$ for all $k\in\ZZ_+$. For
$\y\in\Yq\bez\{-q^{-2k}:\:k\in\ZZ_+\}$ we put
\[
\Fq(\y)=\prod_{k=0}^{\infty}\frac{1+\overline{q^{2k}\y}}{1+q^{2k}\y}.
\]
The infinite product is convergent since $|q|<1$.

We have the following simple 

\begin{Prop}
The function $\Fq$ extends to a continuous function $\Ybar\to\TT$. 
With this extension we have $\Fq(0)=1$ and
\begin{equation}\label{wzorek}
(1+\y)\Fq(\y)=(1+\overline{\y})\Fq(q^2\y)
\end{equation}
for all $\y\in\Ybar$.
\end{Prop}

\begin{Rem}
Dividing both sides of \refeq{wzorek} by $(1+\y)$ and calculating the limit
$\y\to-1$, one finds that
\[
\Fq(-1)=\frac{-\rh}{\rh}\Fq(-q^2).
\]  
This formula will prove useful.
\end{Rem}

The function $\Fq$ was first introduced in \cite{OpEq} with a real deformation
parameter $q$. Its remarkable properties have been very useful in developing
examples of quantum groups (c.f.~e.g.~\cite{e2}). 
The name {\em quantum exponential 
function}\/ is taken from \cite{qef} and will be justified in Subsection
\ref{SumSubs}.

For any function $f$ defined on $\Ybar$ and a fixed $\y\in\Ybar$ we can consider 
a function $\An{f}{\y}$ of a real variable given by
\[
\an{f}{\y}{t}=f(q^{it}\y).
\] 
Throughout the paper much attention will be given to analytic continuation of
such functions.

\begin{Thm}\label{glowne1}
\noindent\begin{itemize}
\item[{\rm (1)}]
For any $\y\in\Ybar$ the function $\An{\Fq}{\y}$ has a holomorphic
continuation to the lower half plane. 
\item[{\rm (2)}] We have the following formula
\begin{equation}\label{rownFq}
\an{\Fq}{q\y}{-i}=(1+\y)\Fq(\y).
\end{equation}
\item[{\rm (3)}] For any $\y\in\Yq$, $\eps>0$ and $M>0$ there exists an $R>0$ 
such that for all $\tau\in\RR$ with $|\tau|<M$ and all $t>R$ 
\[
\bigl|\an{\Fq}{\y}{-t-i\tau}-1\bigr|<\eps.
\]
\end{itemize}
\end{Thm}

\proof
{\sc Ad (1).}
Let $\y=q^nq^{ir}$. Denote
\[
\begin{array}{r@{\;=\;}l@{\smallskip}}
A(t)&\displaystyle{
\Prod_{k=0}^\infty\left(1+\Exp{\frac{2k+n-i(t+r)}{\rhbar}}\right),}\\
B(t)&\displaystyle{
\Prod_{k=0}^\infty\left(1+\Exp{\frac{2k+n+i(t+r)}{\rh}}\right),}
\end{array}
\]
with $t\in\CC$. We can now write $\an{\Fq}{\y}{t}$ as a quotient
\begin{equation}\label{iloraz}
\an{\Fq}{\y}{t}=\frac{A(t)}{B(t)}.
\end{equation}
The functions $A$\/ and $B$\/ extend holomorphically o the whole complex
plane. All zeros of both these functions are simple.

Zeros of $A$\/ form the set
\begin{equation}\label{zeraA}
\left(2\ZZ+1\right)\pi\rhbar-i\left(2\ZZ_++n\right)-r,
\end{equation}
while zeros of $B$\/ form the set
\begin{equation}\label{zeraB}
\left(2\ZZ+1\right)\pi\rh+i\left(2\ZZ_++n\right)-r.
\end{equation}
Let $z=(2p+1)\pi\rh+i(2l+n)-r$ be a zero of $B$ such that 
$\im{z}\leq0$. Notice that
\[
\rh=\rhbar+\rh-\rhbar=\rhbar+2i\iim{\rh}=\rhbar+i\ulamek{N}{\pi}.
\]
Therefore
\[
\begin{array}{r@{\;=\;}l@{\smallskip}}
z&(2p+1)\pi\rh+i(2l+n)-r\\
&(2p+1)\pi\rhbar+i(2p+1)N+i(2l+n)-r\\
&(2p+1)\pi\rhbar+i\bigl((2p+1)N+2l+n\bigr)-r\\
&(2p+1)\pi\rhbar+i(2pN+N+2l+2n-n)-r\\
&(2p+1)\pi\rhbar+i\left(2\left[pN+\ulamek{N}{2}+l+n\right]-n\right)-r\\
&(2p+1)\pi\rhbar-i\left(2\left[-pN-\ulamek{N}{2}-l-n\right]+n\right)-r.
\end{array}
\]
The number $-\left[pN+\ulamek{N}{2}+l+n\right]$ is a positive integer, since
\[
\begin{array}{r@{\;=\;}l@{\smallskip}}
0\geq\im{z}&(2p+1)\pi\iim{\rh}+2l+n\\
&(2p+1)\pi\frac{N}{2\pi}+2l+n\\
&pN+\frac{N}{2}+2l+n\\
&\left[pN+\frac{N}{2}+l+n\right]+l
\end{array}
\]
shows that
\[
-\left[pN+\ulamek{N}{2}+l+n\right]\geq l\in\ZZ_+.
\]
This means that $z$ is a zero of the function $A$. In particular all
singularities of $\An{\Fq}{\y}$ in the lower half plane are removable.
Consequently $\An{\Fq}{\y}$ extends holomorphically to the lower half plane. 

{\sc Ad (2).} To prove formula \refeq{rownFq} denote
\[
\begin{array}{r@{\;=\;}l@{\smallskip}}
\widetilde{A}(t)&\displaystyle{
\Prod_{k=0}^\infty\left(1+\Exp{\frac{2k+n+1-i(t+r)}{\rhbar}}\right),}\\
\widetilde{B}(t)&\displaystyle{
\Prod_{k=0}^\infty\left(1+\Exp{\frac{2k+n+1+i(t+r)}{\rh}}\right)}
\end{array}
\]
Then
\[
\an{\Fq}{q\y}{t}=\frac{\widetilde{A}(t)}{\widetilde{B}(t)}.
\]
Moreover
\[
\begin{array}{r@{\;=\;}l@{\smallskip}}
\widetilde{A}(-i)&\displaystyle{
\Prod_{k=0}^\infty\left(1+\Exp{\frac{2k+n-ir}{\rhbar}}\right)
=\Prod_{k=0}^\infty\left(1+\overline{q^{2k}\y}\right),}\\
\widetilde{B}(-i)&\displaystyle{
\Prod_{k=0}^\infty\left(1+\Exp{\frac{2(k+1)+n+ir}{\rh}}\right)
=\Prod_{k=1}^\infty\left(1+q^{2k}\y\right),}
\end{array}
\]
which shows that
\[
\an{\Fq}{q\y}{-i}=\frac
{\Prod_{k=0}^\infty\left(1+\overline{q^{2k}\y}\right)}
{\Prod_{k=1}^\infty\left(1+q^{2k}\y\right)}
=(1+\y)\prod_{k=0}^\infty\frac{1+\overline{q^{2k}\y}}{1+q^{2k}\y}
=(1+\y)\Fq(\y).
\]

{\sc Ad (3).}
The function 
\[
w\longmapsto
\frac{\prod\limits_{k=0}^\infty\left(1+\overline{q^{2k}w}\right)}
{\prod\limits_{k=0}^\infty\left(1+q^{2k}w\right)}
=\PHase{\prod_{k=0}^{\infty}(1+q^{2k}w)}^{-2}
\]
is continuous in a neighbourhood of $0$ and its value converges to 1 as $w$
tends to 0. Now if $w=\y q^{i(t+i\tau)}$ with $|\tau|$ bounded by $M$ then $w$ 
tends to $0$ as $t$ goes to $-\infty$ and the result follows.  
\qed

\begin{Rem}
The function $\Fq$ can be defined for any non zero $q$ of absolute value
strictly less than $1$. Then for general $\y\in\Yq$ the function
$\An{\Fq}{\y}$ has a meromorphic continuation to the lower half plane. This
continuation is holomorphic if and only if the imaginary part of the inverse
of the logarithm of $q$ is $\ulamek{2\pi}{N}$ with $N\in2\ZZ\bez\{0\}$
(\cite[Twierdzenie 5.11.1]{phd}).
\end{Rem}

\subsection{The product formula}

We shall devote this subsection to proving the following theorem:

\begin{Thm}\label{FFALT}
For any $\y\in\Yq$ we have
\begin{equation}\label{FFal}
\Fq(\y)\Fq(q^2\y^{-1})=\Cq\alpha(q^{-1}\y),
\end{equation}
where
\[
\Cq=e^{-\frac{i}{2}\iim{\rhinv}}\Fq(1)^2.
\]
\end{Thm}

We shall examine the function
\begin{equation}\label{funkcja}
\Yq\ni\y\longmapsto\Fq(\y)\Fq(q^2\y^{-1}).
\end{equation}
Writing $\y=q^nq^{it}$ we shall treat the right hand side of \refeq{funkcja}
as a meromorphic function of $t\in\CC$. It can be rewritten as
\begin{equation}\label{phi1}
\CC\ni t\longmapsto\prod_{k=0}^\infty
\frac{1+\qbar^{2k}\qbar^{n-it}}{1+q^{2k}q^{n+it}}
\prod_{k=1}^\infty\frac{1+\qbar^{2k}\qbar^{-n+it}}{1+q^{2k}q^{-n-it}}
=\frac{\overline{\ph_n(\overline{t})}}{\ph_n(t)},
\end{equation}
where
\[
\ph_n(t)=\prod_{k=0}^\infty\left(1+q^{2k}q^{n+it}\right)\prod_{k=1}^\infty
\left(1+q^{2k}q^{-n-it}\right).
\]
The zeros of the entire function $\ph_n$ are all simple and constitute the set
\[
\Lambda=2\pi\rh\left(\ZZ+\ulamek{1}{2}\right)+i(2\ZZ+n).
\]
It is easy to see that thanks to the form of the imaginary part of $\rh$ the
set $\Lambda$ satisfies $\Lambda=\overline{\Lambda}$. In particular the
function
\[
\CC\ni t\longmapsto\frac{\overline{\ph_n(\overline{t})}}{\ph_n(t)}
\]
is entire.

\begin{Rem}
In fact we have 
\[
\Bigl(\Lambda=\overline{\Lambda}\Bigr)\Longleftrightarrow
\left(\begin{array}{c}
\im{\rh}=\frac{N}{2\pi},\text{ with}\\ N\in2\ZZ
\end{array}\right)
\]
(cf.~\cite[Lemat A.1]{phd}).
\end{Rem}

By the Weierstrass theorem we have the following expression for $\ph_n$:
\begin{equation}\label{postacphi1}
\ph_n(t)=G_n(t)\prod_{\lambda\in\Lambda}\left(1-\ulamek{t}{\lambda}
\right)
e^{\frac{t}{\lambda}+\frac{t^2}{2\lambda^2}},
\end{equation}
where the infinite product is absolutely convergent and $G_n$ is a nowhere
vanishing entire function.

On the other hand it is easy to check that
\[
\begin{array}{r@{\;=\;}l@{\smallskip}}
\ph_n(t)
&\prod\limits_{k\in\ZZ}\left(
1+\exp{\left(\frac{2k+n+it}{\rh}\right)}\right)
\times\left\{\begin{array}{l@{\quad\mathrm{for}\quad}l}
1&k\geq0\\
\exp{\left(-\frac{2k+n+it}{\rh}\right)}&k<0
\end{array}\right.\\
&\prod\limits_{k\in\ZZ}\exp{\left(\frac{2k+n+it}{2\rh}\right)}
2\cosh{\left(\frac{2k+n+it}{2\rh}\right)}
\times\left\{\begin{array}{l@{\quad\mathrm{for}\quad}l}
1&k\geq0\\
\exp{\left(-\frac{2k+n+it}{\rh}\right)}&k<0
\end{array}.\right.
\end{array}
\]
So introducing the function
\[
s(k)=\left\{\begin{array}{l@{\quad\mathrm{for}\quad}l}
1&k\geq0\\
-1&k<0
\end{array}\right.
\]
we obtain
\[
\begin{array}{r@{\;=\;}l@{\smallskip}}
\ph_n(t)
&\prod\limits_{k\in\ZZ}\exp{\left(s(k)\frac{2k+n+it}{2\rh}\right)}
2\cosh{\left(\frac{2k+n+it}{2\rh}\right)}\\
&\prod\limits_{k\in\ZZ}\exp{\left(s(k)\frac{2k+n+it}{2\rh}\right)}
2\left(\prod\limits_{p\in\ZZ}\left(1-\frac{2k+n+it}
{2\pi i\rh(p+\frac{1}{2})}\right)
\exp{\left(\frac{2k+n+it}{2\pi i\rh(p+\frac{1}{2})}\right)}\right).
\end{array}
\]
This can be rewritten as
\begin{equation}\label{postacphi2}
\begin{array}{r@{\;}l}
\ph_n(t)=&\prod\limits_{k\in\ZZ}2\exp{\left(s(k)\frac{2k+n+it}
{2\rh}\right)}\\
&\times\left(\prod\limits_{\lambda\in\Lambda_k}
\left(1-\frac{t}{\lambda}\right)
\left(1+\frac{i(2k+n)}{\lambda-i(2k +n)}\right)
\exp{\left(-\frac{i(2k+n)}{\lambda-i(2k+n)}\right)}
\exp{\left(\frac{t}{\lambda-i(2k+n)}\right)}\right),
\end{array}
\end{equation}
where
\[
\Lambda_k=2\pi\rh\left(\ZZ+\ulamek{1}{2}\right)+i(2k+n).
\]
Of course we have $\Lambda=\bigsqcup\limits_{k\in\ZZ}\Lambda_k$. 
Taking into account \refeq{postacphi1} and \refeq{postacphi2} we obtain
\[
\begin{array}{r@{\;}l}
G_n(t)=&\prod\limits_{k\in\ZZ}2\exp{\left(s(k)\frac{2k+n+it}{2\rh}\right)}\\
&\times\left(\prod\limits_{\lambda\in\Lambda_k}
\left(1+\frac{i(2k+n)}{\lambda-i(2k+n)}\right)
\exp{\left(-\frac{i(2k+n)}{\lambda-i(2k+n)}\right)}
\exp{\left(\frac{it(2k+n)}{\lambda(\lambda-i(2k+n))}-
\frac{t^2}{2\lambda^2}\right)}\right).
\end{array}
\]

For computational reasons it will be more convenient to work with
\[
\begin{array}{r@{\;=\;}l@{\smallskip}}
\frac{G_n(t)}{G_n(0)}
&\prod\limits_{k\in\ZZ}\exp{\left( s(k)\frac{it}{2\rh}\right)}
\left(\prod\limits_{\lambda\in\Lambda_k}
\exp{\left(-\frac{t^2}{2\lambda^2}\right)}
\prod\limits_{\lambda\in\Lambda_k}
\exp{\left(\frac{it(2k+n)}{\lambda(\lambda-i(2k+n))}\right)}\right)\\
&\prod\limits_{k\in\ZZ}\exp{\left(s(k)\frac{it}{2\rh}\right)}
\exp{\left(-\frac{t^2}{2}\Sum_{\lambda\in\Lambda_k}
\frac{1}{\lambda^2}\right)}
\exp{\left(i(2k+n)t\Sum_{\lambda\in\Lambda_k}
\frac{1}{\lambda(\lambda-i(2k+n))}\right)},
\end{array}
\]
especially since from \refeq{postacphi1} we immediately get an expression for
the constant $G_n(0)$: 
\begin{equation}\label{Gzero}
G_n(0)=\prod_{k=0}^\infty(1+q^{2k+n})\prod_{k=1}^\infty(1+q^{2k-n}).
\end{equation}
Moreover using standard methods of complex analysis one can compute the sums
\[
\begin{array}{r@{\;=\;}l@{\smallskip}}
\Sum_{\lambda\in\Lambda_k}\frac{1}{\lambda(\lambda-i(2k+n))}
&\frac{1}{(2k+n)}\frac{1}{2\rh}\tanh{\left(\frac{2k+n}{2\rh}\right)},\\
\Sum_{\lambda\in\Lambda_k}\frac{1}{\lambda^2}
&\left(\frac{1}{2\rh}\right)^2\frac{1}{\left(\cosh{\left(\frac{2k+n}
{2\rh}\right)}\right)^2},
\end{array}
\]
which show that
\[
\ulamek{G_n(t)}{G_n(0)}
=\exp{\left(
\ulamek{it}{2\rh}
\sum_{k\in\ZZ}
\left(
s(k)+\tanh{\left(\ulamek{2k+n}{2\rh}\right)}
\right)
\right)}
\exp{\left(
-\ulamek{t^2}{2}
\left(\ulamek{1}{2\rh}\right)^2
\sum_{k\in\ZZ}\frac{1}
{\left(
\cosh{\left(\frac{2k+n}{2\rh}\right)}
\right)^2}
\right)}.
\]

\begin{Lem}\label{nm1}
We have
\[
\sum_{k\in\ZZ}
\left(
s(k)+\tanh{\left(\frac{2k+n}{2\rh}\right)}
\right)=-(n-1).
\]
\end{Lem}

\proof
First of all notice that since $\re{\rh}<0$ the series in question is
absolutely convergent. To compute its sum we shall use the following,
obvious, formulae:
\begin{eqnarray}
\sum_{k\in\ZZ}\left(s(k)+\tanh{\left(\ulamek{k}{\rh}\right)}
\right)&=&1,\label{C}\\
\sum_{k\in\ZZ+\frac{1}{2}}\left(s(k+\ulamek{1}{2})
+\tanh{\left(\ulamek{k}{\rh}\right)}
\right)&=&0.\label{B}
\end{eqnarray}

Let us deal with the cases of even and odd $n$ separately. First let $n=2l$.
Then
\begin{equation}\label{Aprim}
\sum_{k\in\ZZ}\left(s(k)+\tanh{\left(\ulamek{2k+n}{2\rh}\right)}\right)
=\sum_{k\in\ZZ}\left( s(k-l)+\tanh{\left(\ulamek{k}{\rh}\right)}\right).
\end{equation}
Subtracting the left hand side of \refeq{C} from the right hand side of
\refeq{Aprim} we obtain
\[
\begin{array}{r@{\;}c@{\;}l}
\Sum_{k\in\ZZ}\left(s(k)+\tanh{\left(\frac{2k+n}{2\rh}\right)}\right)
&=&\Sum_{k\in\ZZ}\left(s(k-l)+\tanh{\left(\frac{k}{\rh}\right)}\right)\\
&&\quad-\Sum_{k\in\ZZ}\left(s(k)+\tanh{\left(\frac{k}{\rh}\right)}
\right)+1\\
&=&\Sum_{k\in\ZZ}\left(s(k-l)-s(k)\right)+1=-2l+1\\
&=&-(2l-1)=-(n-1).
\end{array}
\]

Assume now that $n=2l+1$. Then
\begin{equation}\label{Abis}
\sum_{k\in\ZZ}\left(s(k)+\tanh{\left(\ulamek{2k+n}{2\rh}\right)}\right)
=\sum_{k\in\ZZ+\frac{1}{2}}\left(s\left(k-l-\ulamek{1}{2}\right)+\tanh{\left(
\ulamek{k}{\rh}\right)}\right).
\end{equation}
As before we subtract the left hand side of \refeq{B} from the right hand side
of \refeq{Abis} and arrive at
\[
\begin{array}{r@{\;}c@{\;}l}
\Sum_{k\in\ZZ}\left(s(k)+\tanh{\left(\frac{2k+n}{2\rh}\right)}\right)
&=&\Sum_{k\in\ZZ+\frac{1}{2}}\left(s\left(
k-l-\ulamek{1}{2}\right)+\tanh{\left(\frac{k}{\rh}\right)}\right)\\
&&\quad-\Sum_{k\in\ZZ+\frac{1}{2}}\left(
s\left(k-\frac{1}{2}\right)+\tanh{\left(\frac{k}{\rh}\right)}\right)\\
&=&\Sum_{k\in\ZZ+\frac{1}{2}}\left(s\left( k-l-\frac{1}{2}\right)
-s\left(k-\frac{1}{2}\right)\right)\\
&=&\Sum_{k\in\ZZ}\left(s(k-l)-s(k)\right)=-2l=-(n-1),
\end{array}
\]
which ends the proof of Lemma \ref{nm1}.
\qed

Let us introduce the notation
\[
\Theta_n=\sum_{k\in\ZZ}\frac{1}{\left(\cosh{\left(\frac{2k+n}{2\rh}\right)}
\right)^2}.
\]
We can summarize the analysis we have done so far in the following way:
\begin{equation}\label{GnG}
\frac{G_n(t)}{G_n(0)}=\exp{\left(-(n-1)\frac{it}{2\rh}\right)}
\exp{\left(-\frac{t^2}{2}\left(\frac{1}{2\rh}\right)^2\Theta_n\right)}.
\end{equation}

Now let us return to the study of the function \refeq{funkcja}. 
Using \refeq{phi1} and \refeq{postacphi1} we obtain
\[
\Fq(q^nq^{it})\Fq(q^{2-n}q^{-it})=\frac{\overline{\ph_n(\overline{t})}}
{\ph_n(t)}
=\frac{\overline{G_n(\overline{t})}
\prod\limits_{\lambda\in\overline{\Lambda}}\left( 1-\frac{t}{\lambda}\right)
\exp{\left(\frac{t}{\lambda}+\frac{t^2}{2\lambda^2}\right)}}
{G_n(t)\prod\limits_{\lambda\in\Lambda}\left( 1-\frac{t}{\lambda}\right)
\exp{\left(\frac{t}{\lambda}+\frac{t^2}{2\lambda^2}\right)}}
=\frac{\overline{G_n(\overline{t})}}{G_n(t)}.
\]
By \refeq{GnG} this means that
\begin{equation}\label{funkcjastale}
\begin{array}{r@{\;}c@{\;}l}
\Fq(q^nq^{it})\Fq(q^{2-n}q^{-it})&=&
\Phase{G_n(0)}^{-2}\\
&&\times\exp{\left(i(n-1)\Rre{\frac{1}{\rh}}t+
i\Iim{\frac{\Theta_n}{(2\rh)^2}}t^2\right)}.
\end{array}
\end{equation}
In particular the left hand side of \refeq{funkcjastale} extends to an entire
function of $t$. Let us denote this extension by $t\mapsto\Phi(n,t)$.

\begin{Lem}\label{sztuczka}
For any $t\in\RR$  we have
\begin{equation}\label{wzorPhi}
\Phi(n+1,t-i)=q^nq^{it}\Phi(n,t).
\end{equation}
\end{Lem}

\proof
Just as in the proof of Statement (2) of Theorem \ref{glowne1} one can show
that for $\y'\in\Ybar$
\begin{equation}\label{rownFq2}
\an{\Fq}{q\y'}{i}=(1+\y')^{-1}\Fq(q^2\y').
\end{equation}
With $\y=q^nq^{it}$ we have
\[
\Phi(n,t)=\Fq(\y)\Fq(q^2\y^{-1})
\]
and
\[
\Phi(n+1,t+s)=\Fq(q\y q^{is})\Fq(q\y^{-1}q^{-is})
=\an{\Fq}{q\y}{s}\an{\Fq}{q\y^{-1}}{-s}.
\]
Now using \refeq{rownFq} and \refeq{rownFq2} we get
\[
\Phi(n+1,t-i)=\an{\Fq}{q\y}{-i}\an{\Fq}{q\y^{-1}}{i}
=\ulamek{1+\y}{1+\y^{-1}}\Fq(\y)\Fq(q^2\y^{-1})=\y\Phi(n,t)
\]
for $\y\neq-1$. For $\y=-1$ we use the continuity of both sides with respect
to $\y$.
\qed 

Lemma \ref{sztuczka} allows us to determine the constants appearing in 
\refeq{funkcjastale} with simple recurrence formulae. From now on let $t$ 
be a real number. By \refeq{funkcjastale} the right hand side of 
\refeq{wzorPhi} is equal to
\[
\begin{array}{r@{\;}c@{\;}l}
\mathrm{RHS}
&=&q^nq^{it}\Phase{G_n(0)}^{-2}\exp{\left(
i(n-1)\Rre{\frac{1}{\rh}}t+i\Iim{\frac{\Theta_n}{(2\rh)^2}}t^2\right)}\\
&=&\Phase{G_n(0)}^{-2}\exp{\left(
n\Rre{\frac{1}{\rh}}+in\Iim{\frac{1}{\rh}}-\Im{\frac{1}{\rh}}t+i\,
\Re{\frac{1}{\rh}}t
\right)}\\
&&\times\exp{\left(
i(n-1)\Rre{\frac{1}{\rh}}t+i\Iim{\frac{\Theta_n}{(2\rh)^2}}t^2
\right)},
\end{array}
\]
while the left hand side equals
\[
\begin{array}{r@{\;}c@{\;}l}
\mathrm{LHS}
&=&\Phase{G_{n+1}(0)}^{-2}\exp{\left(
in\Rre{\frac{1}{\rh}}(t-i)+
i\Iim{\frac{\Theta_{n+1}}{(2\rh)^2}}(t-i)^2\right)}\\
&=&\Phase{G_{n+1}(0)}^{-2}\exp{\left(n\Rre{\frac{1}{\rh}}+
in\Rre{\frac{1}{\rh}}t-i\Iim{\frac{\Theta_{n+1}}{(2\rh)^2}}\right)}\\
&&\times\exp{\left(2\Iim{\frac{\Theta_{n+1}}{(2\rh)^2}}t+
i\Iim{\frac{\Theta_{n+1}}{(2\rh)^2}}t^2\right)}.
\end{array}
\]
Comparing these expression gives
\[
\begin{array}{l@{\smallskip}}
\Phase{G_{n}(0)}^{-2}\exp{\left(in\Iim{\frac{1}{\rh}}-\Im{\frac{1}{\rh}}t+
i\Iim{\frac{\Theta_n}{(2\rh)^2}}t^2\right)}\\
\qquad=\Phase{G_{n+1}(0)}^{-2}\exp{\left(-i\,
\Im{\frac{\Theta_{n+1}}{(2\rh)^2}}+2\Iim{\frac{\Theta_{n+1}}{(2\rh)^2}}t
+i\Iim{\frac{\Theta_{n+1}}{(2\rh)^2}}t^2
\right)}
\end{array}
\]
for all $t\in\RR$. As first and second derivatives of
$\mathrm{LHS}$ and $\mathrm{RHS}$ at $t=0$ are equal, we obtain
\begin{equation}\label{rownosci}
\Im{\frac{\Theta_n}{(2\rh)^2}}=\Im{\frac{\Theta_{n+1}}{(2\rh)^2}}=
-\frac{1}{2}\Iim{\frac{1}{\rh}}.
\end{equation}
At the same time the equality of values of $\mathrm{LHS}$ and $\mathrm{RHS}$
at $t=0$ shows that
\begin{equation}\label{rownosci2}
\Phase{G_n(0)}^{-2}\exp{\left( in\Iim{\ulamek{1}{\rh}}\right)}
=\Phase{G_{n+1}(0)}^{-2}\exp{
\left(-i\Iim{\ulamek{\Theta_{n+1}}{(2\rh)^2}}\right)}.
\end{equation}

Using \refeq{rownosci} and \refeq{rownosci2} we find the recurrence relation
\[
\Phase{G_{n+1}(0)}^{-2}=\exp{\left(i\left(n-\ulamek{1}{2}\right)
\Iim{\ulamek{1}{\rh}}\right)}\Phase{G_n(0)}^{-2},
\]
which we turn into
\begin{equation}\label{wzornaPhase}
\Phase{G_n(0)}^{-2}=\exp{\left( -\ulamek{i}{2}\Iim{\ulamek{1}{\rh}}\right)}
\exp{\left( i\ulamek{(n-1)^2}{2}\Iim{\ulamek{1}{\rh}}\right)}
\Phase{G_0(0)}^{-2}.
\end{equation}

Now inserting \refeq{rownosci} and \refeq{wzornaPhase} into
\refeq{funkcjastale} we get
\[
\begin{array}{r@{}l}
\Fq(q^nq^{it})&\Fq(q^{2-n}q^{-it})=\exp{\left(-\frac{i}{2}
\Iim{\frac{1}{\rh}}\right)}\\
&\times\Phase{G_0(0)}^{-2}\exp{\left(\frac{i}{2}\left((n-1)^2-t^2\right)
\Iim{\frac{1}{\rh}}+i(n-1)t\Rre{\frac{1}{\rh}}\right)},
\end{array}
\]
which in view of \refeq{Alfa} means that
\[
\Fq(\y)\Fq(q^2\y^{-1})=\exp{\left(-\ulamek{i}{2}\Iim{\ulamek{1}{\rh}}\right)}
\Phase{G_0(0)}^{-2}\alpha(q^{-1}\y),
\]
where $\y=q^nq^{it}$.

Finally equation \refeq{Gzero} together with the definition of $\Fq$
gives the expression for the constant $\Cq$ in \refeq{FFal}. 

This concludes the proof of Theorem \ref{FFALT}. The name {\em product
formula}\/ for Equation \refeq{FFal} is self explanatory. It is an interesting
fact that for other values of the deformation parameter $q$ than those
considered in this paper, the proof of the analogous formula simplifies
considerably (cf.~\cite[Section 1]{OpEq} and \cite[Section 1]{azb}).

\subsection{Asymptotic behaviour of $\Fq$}

\begin{Prop}\label{wzrostFq}
The function $\Fq$ has the following asymptotic behaviour: for any
$\y\in\Ybar$ and $t,\tau\in\RR$ we have
\[
\bigl|\an{\Fq}{\y}{t-i\tau}\bigr|=\Xi_\y(t-i\tau)|q^{-1}q^{it}\y|^{\tau},
\]
where
\[
\lim_{t\to\infty}\Xi_\y(t-i\tau)=1
\]
for any $\tau\in\RR$.
\end{Prop}

\proof
The mapping
\[
t\longmapsto\an{\Fq}{\y}{t}\an{\Fq}{q^{2}\y^{-1}}{-t}=
\Fq(q^{it}\y)\Fq\bigl(q^2(q^{it}\y)\bigr)=\Cq\alpha(q^{-1}q^{it}\y)
\]
extends to an entire function on $\CC$. Statement (3) of Theorem
\ref{glowne1} says that with bounded $\tau$ we have
\[
\an{\Fq}{q^2\y^{-1}}{-t-i\tau}\xrightarrow[t\to\infty]{}1.
\]
Denoting the reciprocal of the absolute value of this function by $\Xi_\y$ we
obtain
\[
\bigl|\an{\Fq}{\y}{t-i\tau}\bigr|\Xi_\y(t-i\tau)^{-1}
=\left|\an{\alpha}{q^{-1}\y}{t-i\tau}\right|.
\]
It remains to determine the asymptotic behaviour of the analytic continuation
of $\alpha$. With $\y=q^mq^{is}$ we have (cf.~\refeq{Alfa})
\[
\begin{array}{r@{\;=\;}l@{\smallskip}}
|\an{\alpha}{q^{-1}\y}{t-i\tau}|
&\left|e^{\frac{i}{2}\bigl((m-1)^2-(s+t-t\tau)^2\bigr)\im{\rh^{-1}}}
e^{i(m-1)(s+t-i\tau)\re{\rh^{-1}}}\right|\\
&e^{-(s+t)\tau\iim{\rh^{-1}}}e^{\tau(m-1)\re{\rh^{-1}}}\\
&|q^{-1}q^mq^{is}q^{it}|^\tau=|q^{-1}\y q^{it}|^\tau.
\end{array}
\]
\qed

\subsection{Uniqueness of $\Fq$}

\begin{Prop}\label{uniq}
Let $\Phi$ be a continuous function on $\Ybar$ such that
\begin{itemize}
\item[{\rm (1)}] $\Phi(0)=1$,
\item[{\rm (2)}] For any $\y\in\Ybar$ the function 
$\An{\Phi}{\y}\colon\RR\ni t\mapsto\Phi(q^{it}\y)\in\CC$ has a holomorphic
extension to the lower half plane,
\item[{\rm (3)}] for any $\y\in\Ybar$ we have
\begin{equation}\label{rownPhi}
\an{\Phi}{q\y}{-i}=(1+\y)\Phi(\y),
\end{equation}
\item[{\rm (4)}] For any $\delta>0$ and any $\y\in\Ybar$ there exist constants 
$C_1$
and $C_2$ such that for any $0\leq\sigma<1$ and any $s\in\RR$ we have
\[
\bigl|\an{\Phi}{\y}{s-i\sigma}\bigr|
\leq C_1+C_2\left|q^{-1}\y q^{is}\right|^{1+\delta}.
\]
\end{itemize}
Then $\Phi=\Fq$
\end{Prop}

\proof
Let us fix a $\y\in\Ybar$ and define
\[
\ph_\y(z)=\frac{\an{\Phi}{\y}{z}}{\an{\Fq}{\y}{z}}
\]
for $z$ in the lower half plane. This way we obtain a meromorphic function
$\ph_\y$ on the lower half plane. Our aim is to prove that it is a constant
function equal to one. We shall show this in four steps:
\begin{enumerate}
\item $\ph_\y$ is holomorphic;
\item $\ph_\y$ is periodic with non real period, in particular $\ph_\y$
extends to an entire function;
\item $\ph_\y$ factorizes through the map $s\mapsto q^{is}$, more precisely
there exists an entire function $\psi_\y$ such that
\[
\ph_\y(s-i\sigma)=\psi_\y\bigl(q^{i(s-i\sigma)}\bigr);
\]
\item $\psi_\y$ is constant equal to $1$.
\end{enumerate}

{\sc Ad 1.} Analysis of zeros of the function $\An{\Fq}{\y}$ 
(cf.~\refeq{iloraz},
\refeq{zeraA} and \refeq{zeraB}) shows that all zeros of this function have
integer imaginary parts. By \refeq{rownFq} we have
\[
\an{\Fq}{\y}{t-i}=(1+q^{it}q^{-1}\y)\Fq(q^{it}q^{-1}\y)
\]
for all $t\in\RR$. Therefore the only possible zero on the line $\RR-i$ exists
when there is a $t_0\in\RR$ such that $q^{it_0}\y=-q$. By \refeq{rownPhi} we
have the same facts for $\Phi$. Therefore the (only possible) zero of 
$\An{\Fq}{\y}$ in the strip $\{z\in\CC:\:-2<\im{z}\leq0\}$ cancels with one of
zeros of $\An{\Phi}{\y}$. It follows that $\ph_\y$ extends to a holomorphic
function in $\{z\in\CC:\:-2<\im{z}<0\}$. Notice further that by \refeq{rownFq}
and \refeq{rownPhi} we have
\begin{equation}\label{Nrazy}
\ph_\y(t-i)=\frac{\an{\Phi}{\y}{t-i}}{\an{\Fq}{\y}{t-i}}=
\frac{(1+q^{it}q^{-1}\y)\an{\Phi}{q^{-1}\y}{t}}
{(1+q^{it}q^{-1}\y)\an{\Fq}{q^{-1}\y}{t}}=\ph_{q^{-1}\y}(t).
\end{equation}
for all $t\in\RR$ and so this equality remains true for $t$ in the lower half
plane. Therefore for any $k\in\NN$ holomorphy of $\ph_\y$ in the strip 
$\{z\in\CC:\:-(k+2)<\im{z}<k\}$ is equivalent to that of $\ph_{q^{-k}\y}$ in 
$\{z\in\CC:\:-2<\im{z}<0\}$. In particular $\ph_\y$ is holomorphic in the
lower half plane.

{\sc Ad 2.} Using \refeq{Nrazy} $N$ times we get 
\[
\ph_\y(t-Ni)=\ph_{q^{-N}\y}(t)=\ph_{q^{it_0}\y}(t),
\]
for some $t_0\in\RR$ (namely such that $q^{it_0}=q^{-N}$) and all $t\in\RR$.
In particular
\[
\ph_\y(t-Ni)=\frac{\an{\Phi}{q^{it_0}\y}{t}}{\an{\Fq}{q^{it_0}\y}{t}}
=\frac{\an{\Phi}{\y}{t+t_0}}{\an{\Fq}{\y}{t+t_0}}=\ph_\y(t+t_0)
\]
or in other words $\ph_\y\bigl(t-(t_0+Ni)\bigr)=\ph_\y(t)$ for all $t\in\RR$. 
By holomorphy of $\ph_\y$ this equality holds for $t$ in the
lower half plane. 

{\sc Ad 3.} for any $s\in\RR$ we have
\[
\ph_\y(s)=\frac{\an{\Phi}{\y}{s}}{\an{\Fq}{\y}{s}}=
\frac{\Phi(q^{is}\y)}{\Fq(q^{is}\y)}.
\]
and by periodicity
\[
\ph_\y(s)=\ph_\y\bigl(s-(t_0+Ni)\bigr)=
\frac{\Phi\left(q^{i\bigl(s-(t_0+Ni)\bigr)}\y\right)}
{\Fq\left(q^{i\bigl(s-(t_0+Ni)\bigr)}\y\right)}.
\]
As $z$ goes along a path from $s$ to $s-(t_0+Ni)$ the variable $q^{iz}\y$
goes along a closed path beginning and ending in $q^{is}\y$. Therefore the
formula
\[
\psi_\y(q^{iz})=\ph_\y(z)
\]
defines a holomorphic function $\psi_\y$ on $\CC\bez\{0\}$. 
By Statement (3) of Theorem \ref{glowne1} we know that $\an{\Fq}{\y}{z}$
converges to $1$ as the real part of $z$ goes to $-\infty$ and the imaginary
part stays bounded. Also by assumption (4) of this proposition
$\an{\Phi}{\y}{z}$ is bounded when $z$ moves to $-\infty$. Therefore the
quotient is bounded. Consequently $\psi_\y(z)$ is bounded as $z\to 0$. It
follows that $\psi_\y$ is entire. 

{\sc Ad 4.} Let us fix a $0<\delta<\ulamek{1}{2}$. 
We know that there are constants $C_1$ and $C_2$ such that 
\[
\bigl|\an{\Phi}{\y}{s-i\sigma}\bigr|\leq C_1+C_2|q^{-1}\y q^{is}|^{1+\delta}
\]
for all $s\in\RR$ and $0\leq\sigma<1$. By Proposition \ref{wzrostFq}, for 
$s\to\infty$ we have
\begin{equation}\label{ogr}
\bigl|\ph_\y(s-i\sigma)\bigr|\leq\frac{C_1}
{\Xi_\y(s-i\sigma)|q^{-1}\y q^{is}|^\sigma}+\frac{C_1}{\Xi_\y(s-i\sigma)}
|q^{-1}\y q^{is}|^{1+\delta-\sigma}.
\end{equation}
Consider now the values of $\psi_\y$ on the curve 
$\Upsilon=\{q^{i(s-i\sigma)}\y:\:s\in\RR\}$. As $s\to\infty$ the corresponding
point of $\Upsilon$ goes to infinity. Now 
\refeq{ogr} shows that for $z\in\Upsilon$ we have
\[
|\psi_\y(z)|=o(|z|^2).
\]
It follows that 
\[
\psi_\y(z)=\psi_\y(0)+\psi_{\y}'(0)z.
\]
However if $\sigma>1-\delta$ then $1+\delta-\sigma<2\delta<1$ and consequently
\[
\bigl|\psi_\y(z)\bigr|=o(|z|).
\]
In particular $\psi_{\y}'(0)=0$ and $\psi_\y(z)=\psi_\y(0)$ for all $z$.
By assumption (1)
\[
\psi_\y(0)=\lim_{\RR\ni t\to-\infty}\ph_\y(t)
=\lim_{\Ybar\ni\y\to0}\frac{\Phi(\y)}{\Fq(\y)}=1.
\]
\qed

In the next subsection we shall exhibit a function satisfying conditions
(1)--(4) of Proposition \ref{uniq}. In particular that will imply that $\Fq$
satisfies these conditions.

\subsection{Fourier transform of $\Fq$}

All results of this subsection can be proved with more or less standard
techniques from the theory of functions of one complex variable and theory of
distributions. Therefore we have decided to present only sketches of proofs.
The details can have been taken care of in \cite[Uzupe\l{}nienie A.2]{phd} 
(cf.~also \cite[Appendix B]{azb}).

As any locally compact group $\Yq$ possesses a Haar measure $d\Mu$. We shall 
chose the following normalization:
\[
\Int_{\Yq}f(\y)\,d\Mu(\y)
=\sum_{k=0}^{N-1}\Int_{-\infty}^{+\infty}f(q^kq^{it})\,dt.
\] 
Apart from integration we shall also use theory of distributions on $\Yq$. To
that end let  us define the Schwartz space $\Sch(\Yq)$ as the space of all
functions $f\colon\Yq\to\CC$ such that the functions
\[
\RR\ni t\longmapsto f(q^kq^{it})
\]
belong to the space $\Sch(\RR)$ (the usual Schwartz space on $\RR$) for 
$k=1,\ldots,N$. This definition is, of course, compatible with the isomorphism
$\Yq\cong\ZZ_N\times\RR$ (i.e. $\Sch(\Yq)\cong\CC^N\otimes\Sch(\RR)$).

We shall consider the following function
\[
\Hq(\y)=\frac{\chi(-q^{-2},\y)\y}{(1-\overline{\y})\Fq(-q^2\y)}
\]
defined for $\y\in\Yq\bez\{1\}$. This function defines a tempered distribution
on $\Yq$ by integration. The pole of the function $\Hq$ has to be rounded. Let
$\ell$ be the oriented contour in $\CC$ coinciding with $\RR$, but rounding the
point $0$ from above. For $f\in\Sch(\Yq)$ we have
\[
\la\Hq,f\ra=\Int_\ell\Hq(q^{it})f(q^{it})\,dt+
\sum_{k=1}^{N-1}\Int_{-\infty}^{+\infty}\Hq(q^kq^{it})f(q^kq^{it})\,dt.
\] 
We let $\Phiq$ be the inverse Fourier transform of the tempered distribution
$\Hq$:
\begin{equation}\label{Phiq1}
\Phiq(\y)=\Int_\ell
\frac{\chi(-q^{-2}\y,q^{it})q^{it}}
{\left(1-\overline{q^{it}}\right)\Fq(-q^2q^{it})}\,dt+
\sum_{k=1}^{N-1}\Int_{-\infty}^{+\infty}
\frac{\chi(-q^{-2}\y,q^kq^{it})q^kq^{it}}
{\left(1-\overline{q^kq^{it}}\right)\Fq(-q^{k+2}q^{it})}\,dt.
\end{equation}

\begin{Prop}\label{Phiq}
\noindent
\begin{itemize}
\item[{\rm (1)}] The integral \refeq{Phiq1} is convergent as an improper
Riemann integral. More precisely there exists the limit
\begin{equation}\label{Riem}
\begin{array}{r@{\;}c@{\;}l@{\smallskip}}
\Phiq(\y)&=&\lim\limits_{R\to\infty}\INt{\ell^R}
\frac{\chi(-q^{-2}\y,q^{it})q^{it}}
{\left(1-\overline{q^{it}}\right)\Fq(-q^2q^{it})}\,dt\\
&&+\lim\limits_{R\to\infty}\SUM{k=1}{N-1}\INT{-\infty}{R}
\frac{\chi(-q^{-2}\y,q^kq^{it})q^kq^{it}}
{\left(1-\overline{q^kq^{it}}\right)\Fq(-q^{k+2}q^{it})}\,dt,
\end{array}
\end{equation}
where $\ell^R$ is the part of $\ell$ starting at $-\infty$ and ending at $R$.
\item[{\rm (2)}] For any $\tau\in]0,1[$ we have
\begin{equation}\label{Riem2}
\begin{array}{r@{\;}c@{\;}l@{\smallskip}}
\Phiq(\y)&=&-2\pi\rhbar\Fq(-q^2)^{-1}\\
&&+\SUM{k=0}{N-1}\left(q\PHAse{-q^{-2}\y}\right)^k
\INt{\RR-i\tau}\frac{|\y|^{iz}|q|^{-2iz}q^{iz}}
{\left(1-\overline{q^kq^{iz}}\right)\an{\Fq}{-q^{k+2}}{z}}\,dz.
\end{array}
\end{equation}
\item[{\rm (3)}] $\Phiq$ extends to a continuous function on $\Ybar$ and 
\[
\lim_{\y\to0}\Phiq(\y)=-2\pi\rhbar\Fq(-q^2)^{-1}.
\]
\item[{\rm (4)}] For any $\y\in\Ybar$ the function 
$\RR\ni t\mapsto\Phiq(q^{it}\y)$ has a holomorphic continuation to the lower
half plane.
\item[{\rm (5)}] For any $\y\in\Yq$ we have
\[
\an{\Phiq}{q\y}{-i}=(1+\y)\Phiq(\y).
\]
\item[{\rm (6)}] For any $\delta>0$ and any $\y\in\Ybar$ there exist constants 
$C_1$ and $C_2$ such that for any $0\leq\sigma<1$ and any $s\in\RR$ we have
\[
\bigl|\an{\Phiq}{\y}{s-i\sigma}\bigr|
\leq C_1+C_2\left|q^{-1}\y q^{is}\right|^{1+\delta}.
\]
\end{itemize}
\end{Prop}

\proof
{\sc Ad (1).} Due to the exponential decay of $|q^{it}|$ for $t\to-\infty$ the
integral over $\ell^R$ is convergent. After elementary manipulations we can
rewrite the right hand side of \refeq{Riem} as
\[
\lim_{R\to\infty}
\sum_{k=0}^{N-1}\Int_{\ell^R}\frac
{\left(q\PHAse{-q^{-2}\y}\right)^k
e^{it\left(\log{|\y|}-\re{\rhinv}\right)}e^{-t\iim{\rhinv}}}
{\bigl(1-e^{-ik\iim{\rhinv}}e^{-it\rre{\rhinv}}
e^{k\rre{\rhinv}}e^{-t\iim{\rhinv}}\bigr)\an{\Fq}{-q^{k+2}}{t}}\,dt.
\]
Then we choose a number $\tau\in]0,1[$ and deform the integration contour in
the following way:
\begin{figure}[h]
\begin{center}
\begin{pspicture}(-7,-2.5)(7,2.5)
\psaxes[arrowscale=2,labels=none,ticks=none]{->}(0,0)(-6.8,-2.3)(6.95,2.45)
\psline[linewidth=2pt,linearc=.1]{c-c}(-6.5,0)(-.4,0)
\psarc[linewidth=2pt](0,0){.4}{0}{180}
\psline[linewidth=2pt,linearc=.1]{c-c}(.4,0)(1.5,0)(1.5,-1.5)(6,-1.5)(6,0)
\psline[linewidth=2pt,arrowscale=1.5]{->}(-4,0)(-3.3,0)
\put(.1,-.35){$0$}
\psline[linewidth=.5pt]{-}(6,.1)(6,-.1)
\put(5.8,0.2){$R$}
\psline[linewidth=.5pt,linestyle=dotted](0,-1.5)(2,-1.5)
\psline[linewidth=.5pt]{-}(-.1,-1.5)(.1,-1.5)
\put(-.8,-1.55){$-i\tau$}
\psline[linewidth=.5pt]{-}(1.5,.1)(1.5,-.1)
\put(1.4,.2){$1$}
\end{pspicture}
\end{center}
\caption{Deformation of the contour $\ell^R$}
\end{figure}
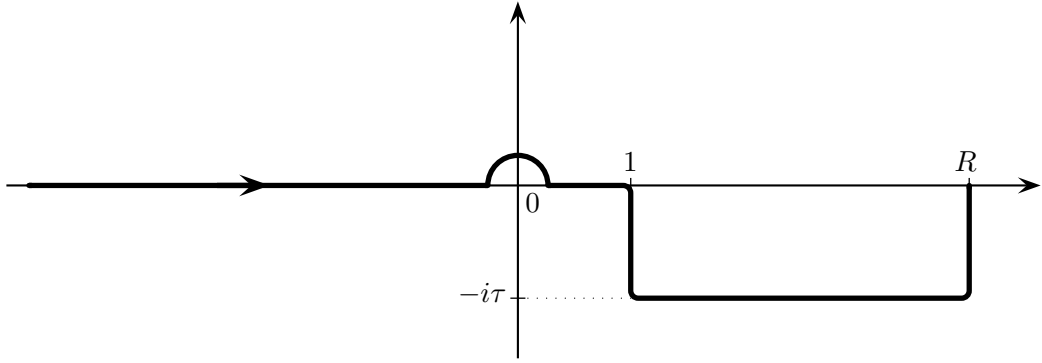
Now standard methods show that the integral along the part of the contour from
$R-i\tau$ to $R$ goes to $0$ as $R\to\infty$ and that the integral over the
remaining part of the contour has a limit as $R\to\infty$.

{\sc Ad (2).} To obtain formula \refeq{Riem2} we deform further the integration
contour:
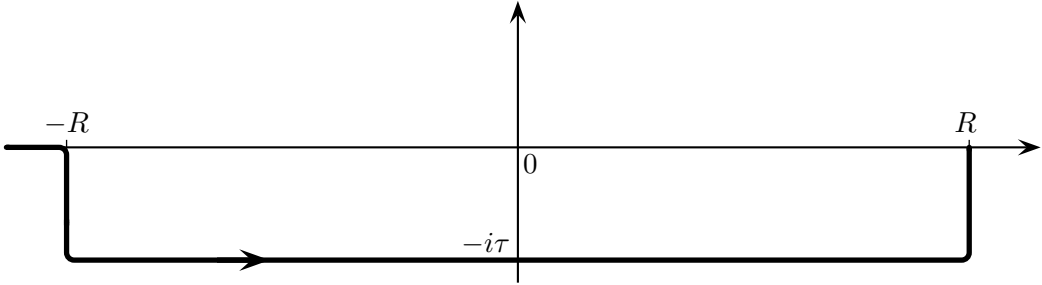
\begin{figure}[H]
\begin{center}
\begin{pspicture}(-7,-2)(7,2)
\psaxes[arrowscale=2,labels=none,ticks=none]{->}(0,0)(-6.8,-1.8)(6.95,1.95)
\psline[linewidth=2pt,linearc=.1]{c-c}(-6.8,0)(-6,0)(-6,-1)
\psline[linewidth=2pt,linearc=.1]{c-c}(-6,-1)(-6,-1.5)(6,-1.5)(6,0)
\psline[linewidth=2pt,arrowscale=1.5]{->}(-4,-1.5)(-3.3,-1.5)
\put(.07,-.35){$0$}
\psline[linewidth=.5pt]{-}(6,.1)(6,-.1)
\put(5.8,0.2){$R$}
\psline[linewidth=.5pt]{-}(-.1,-1.5)(.1,-1.5)
\put(-.75,-1.4){$-i\tau$}
\psline[linewidth=.5pt]{-}(-6,.1)(-6,-.1)
\put(-6.3,.2){$-R$}
\end{pspicture}
\end{center}
\caption{Further deformation of the contour $\ell^R$}
\end{figure}
This deformation does not change the value of the summands for
$k\in\{1,\ldots,N-1\}$. For $k=0$ we obtain $-2\pi\rhbar\Fq(-q^2)^{-1}$ as the
residue of the integrand at the point $0$. Again standard computations show
that The integral over the line from $-R$ to $-R-i\tau$ tends to $0$ as
$R\to\infty$. In the limit we get \refeq{Riem2}.

{\sc Ad (3).} This follows from \refeq{Riem2}.

{\sc Ad (4).} For $\y\in\Ybar$ and $s\in\RR$ the value $\an{\Phiq}{\y}{s}$ is a
limit over $R\to\infty$ of a sum of terms of the form
\begin{equation}\label{holows}
|q^k|^{is}q^k\PHase{-q^{-2}\y}^k\Int_{\ell^R}\frac
{|-q^{-2}\y|^{it}|q^{it}|^{is}q^{it}\,dt}
{\left(1-\overline{q^k q^{it}}\right)\an{\Fq}{-q^{k+2}}{t}}
\end{equation}
Now we want to put $s-i\sigma$ in place of $s$ with $\sigma>0$. To that end we
must deform $\ell^R$ in such a way that the integral be convergent and that we
avoid all zeros of the denominator. An easy analysis of the zeros of
$\An{\Fq}{-q^{k+2}}$ shows that they all lie below the line $\RR\rh+(k+2)i$.
Therefore if we choose $\tau>\sigma$ and deform the integration contour 
as follows:
\begin{figure}[H]
\begin{center}
\begin{pspicture}(-5,-4)(8,3)
\psaxes[arrowscale=2,labels=none,ticks=none]{->}(0,0)(-4.8,-3.8)(7.95,2.95)
\psline[linewidth=1pt](-.5,3)(3,-4) 
\psline[linewidth=2pt]{c-c}(-4.5,0)(-.4,0)
\psarc[linewidth=2pt](0,0){.4}{0}{180}
\psline[linewidth=2pt,linearc=.1]{c-c}(.4,0)(2.2,0)(2.2,-2)(7,-2)(7,0)
\psline[linewidth=2pt,arrowscale=1.5]{->}(-3,0)(-2.3,0)
\rput{-62}(2.85,-3.2){
$\scriptstyle{\RR\rh+(k+2)i}$}
\psline[linewidth=.5pt,linestyle=dotted](0,-2)(2.2,-2)
\psline[linewidth=.5pt]{-}(-.1,-2)(.1,-2)
\put(-.8,-2.05){$-i\tau$}
\psline[linewidth=.5pt]{-}(7,.1)(7,-.1)
\put(6.8,0.2){$R$}
\psline[linewidth=.5pt]{-}(2.2,.1)(2.2,-.1)
\put(2.1,.2){$c$}
\put(2.4,1.7){
$c=-(\tau+k+2)\frac{2\pi\rre{\rh}}{N}+\eps$}
\end{pspicture}
\end{center}
\caption{Avoiding zeros of $\An{\Fq}{-q^{k+2}}$}\label{tutaj}
\end{figure}
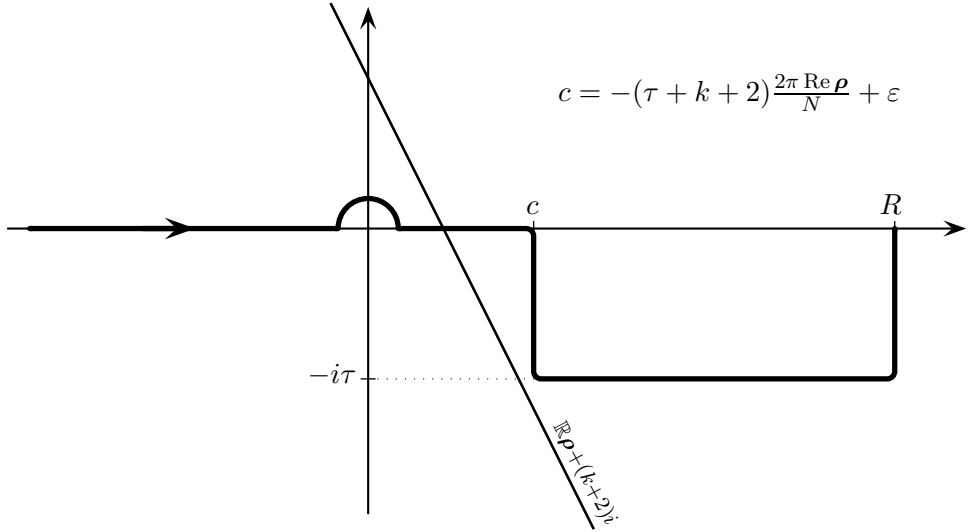
The value of \refeq{holows} will not change. Now, as before, it is possible to
show that the integral over the line from $R-i\tau$ to $R$ goes to $0$ as
$R\to\infty$. At the same time it is easy to see that the limit $R\to\infty$
defines a holomorphic function of $z=s-i\sigma$ for $0<\sigma<\tau$. Since
$\tau$ was arbitrarily large we see that $\An{\Phiq}{\y}$ has a holomorphic
extension to the lower half plane. In particular for $z=-i$ and with $q\y$
in place of $\y$ we have
\begin{equation}\label{przedl}
\an{\Phiq}{q\y}{-i}=\Sum_{k=0}^{N-1}\Int_\ell
\frac{\chi(-q^{-2}\y,q^kq^{it})\left(q^kq^{it}\right)^2\,dt}
{\left(1-\overline{q^kq^{it}}\right)\Fq(-q^{k+2}q^{it})}.
\end{equation}

{\sc Ad (5).} Combining \refeq{wzorek} and \refeq{rownFq} we get 
\[
\an{\Fq}{\y}{-i}=(1-\overline{\y})\Fq(-q^2\y)
=(1-\overline{\y})\an{\Fq}{-q^2\y}{0}
\]
for all $\y\in\Ybar$ and more generally
\[
\an{\Fq}{\y}{t-i}=(1-\overline{q^{it}\y})\an{\Fq}{-q^2\y}{t}
\]
for all $t\in\RR$.

Using this formula we compute
\[
\begin{array}{r@{\;}c@{\;}l@{\smallskip}}
\Phiq(\y)&=&\displaystyle{\Sum_{k=0}^{N-1}\Int_\ell
\frac{\chi(-q^{-2}\y,q^kq^{it})q^kq^{it}\,dt}
{\left(1-\overline{q^kq^{it}}\right)\Fq(-q^{k+2}q^{it})}}\\
&=&\displaystyle{\Sum_{k=0}^{N-1}\Int_\ell
\frac{\PHase{-q^{-2}\y}^k|-q^{-2}\y|^{it}q^kq^{it}\,dt}
{\left(1-\overline{q^kq^{it}}\right)\an{\Fq}{-q^{k+2}}{t}}}\\
&=&\displaystyle{\Sum_{k=0}^{N-1}\Int_\ell
\frac{\PHase{-q^{-2}\y}^k|-q^{-2}\y|^{it}q^kq^{it}\,dt}
{\an{\Fq}{-q^{k+1}}{t-i}}}\\
&=&\displaystyle{\Sum_{k=0}^{N-1}\Int_{\ell-i}
\frac{\PHase{-q^{-2}\y}^k|-q^{-2}\y|^{i(z+i)}q^kq^{i(z+i)}\,dz}
{\an{\Fq}{-q^{k+1}}{z}},}
\end{array}
\]
and so
\[
\begin{array}{r@{\;}c@{\;}l@{\smallskip}}
\Phiq(\y)&=&\displaystyle{\PHase{-q^{-2}\y}^{-1}|-q^{-2}\y|^{-1}}\\
&&\displaystyle{\times
\Sum_{k=0}^{N-1}\Int_{\ell-i}
\frac{\PHase{-q^{-2}\y}^{k+1}
|-q^{-2}\y|^{iz}q^{k-1}q^{iz}\,dz}{\an{\Fq}{-q^{k+1}}{z}}}\\
&=&\displaystyle{\left(-q^{-2}\y\right)^{-1}\Sum_{k=0}^{N-1}\Int_{\ell-i}
\frac{\PHase{-q^{-2}\y}^{k+1}|-q^{-2}\y|^{iz}q^{k-1}q^{iz}\,dz}
{\an{\Fq}{-q^{k+1}}{z}}.}
\end{array}
\]
Now we change the integration contour from $\ell-i$ to $\ell$ to obtain
\[
\begin{array}{r@{\;}c@{\;}l@{\smallskip}}
\Phiq(\y)&=
&\displaystyle{\left(-q^{-2}\y\right)^{-1}\Sum_{k=0}^{N-1}\Int_\ell
\frac{\PHase{-q^{-2}\y}^{k+1}|-q^{-2}\y|^{it}q^{k-1}q^{it}\,dt}
{\an{\Fq}{-q^{k+1}}{t}}}\\
&=&\displaystyle{\left(-q^{-2}\y\right)^{-1}\Sum_{k=0}^{N-1}\Int_\ell
\frac{\chi(-q^{-2}\y,q^kq^{it})q^{k-1}q^{it}\,dt}
{\an{\Fq}{-q^{k+1}}{t}}}\\
&=&\displaystyle{-\y^{-1}\Sum_{k=0}^{N-1}\Int_\ell
\frac{\chi(-q^{-2}\y,q^kq^{it})q^{k+1}q^{it}\,dt}
{\an{\Fq}{-q^{k+1}}{t}}}\\
&=&\displaystyle{-\y^{-1}\Sum_{k=0}^{N-1}\Int_\ell
\frac{\chi(-q^{-2}\y,q^{k+1}q^{it})q^{k+1}q^{it}\,dt}
{\Fq(-q^{k+1}q^{it})}.}
\end{array}
\]
Now by \refeq{wzorek} we have
\[
\begin{array}{r@{\;=\;}l@{\smallskip}}
\Phiq(\y)
&\displaystyle{-\y^{-1}\Sum_{k=0}^{N-1}\Int_\ell
\frac{\chi(-q^{-2}\y,q^{k+1}q^{it})q^{k+1}q^{it}\,dt}
{\frac{1-\overline{q^{k+1}q^{it}}}{1-q^{k+1}q^{it}}
\Fq(-q^2q^{k+1}q^{it})}}\\
&\displaystyle{-\y^{-1}\Sum_{k=0}^{N-1}\Int_\ell
\frac{\chi(-q^{-2}\y,q^{k+1}q^{it})(1-q^{k+1}q^{it})q^{k+1}q^{it}\,dt}
{\left(1-\overline{q^{k+1}q^{it}}\right)\Fq(-q^2q^{k+1}q^{it})}}\\
&\displaystyle{-\y^{-1}\Sum_{k=0}^{N-1}\Int_\ell
\frac{\chi(-q^{-2}\y,q^kq^{it})(1-q^kq^{it})q^kq^{it}\,dt}
{\left(1-\overline{q^kq^{it}}\right)\Fq(-q^{k+2}q^{it})}}\\
&\displaystyle{-\y^{-1}\Sum_{k=0}^{N-1}\Int_\ell
\frac{\chi(-q^{-2}\y,q^kq^{it})q^kq^{it}\,dt}
{\bigl(1-\overline{q^kq^{it}}\bigr)\Fq(-q^{k+2}q^{it})}
+\y^{-1}\Sum_{k=0}^{N-1}\Int_\ell
\frac{\chi(-q^{-2}\y,q^kq^{it})\left(q^kq^{it}\right)^2\,dt}
{\bigl(1-\overline{q^kq^{it}}\bigr)\Fq(-q^{k+2}q^{it})}.}
\end{array}
\]
In other words (cf.~\refeq{przedl})
\[
\an{\Phiq}{\y}{0}=\y^{-1}\an{\Phiq}{q\y}{-i}-\y^{-1}\an{\Phiq}{\y}{0},
\]
i.e.
\[
\an{\Phiq}{q\y}{-i}=(1+\y)\an{\Phiq}{\y}{0}=(1+\y)\Phiq(\y).
\]

{\sc Ad (6).}
$\an{\Phiq}{\y}{s-i\sigma}$ is a sum of terms of the form
\begin{equation}\label{modul}
\left(|q|^{is}q^kq^{\sigma}\PHase{-q^{-2}\y}\right)^k
\Int_{\ell}\frac{|q^{it}|^{is}\,|q^{it}|^{\sigma}\,|-q^{-2}\y|^{it}q^{it}\,dt}
{\left(1-\overline{q^kq^{it}}\right)\an{\Fq}{-q^{k+2}}{t}},
\end{equation}
where the integral is understood as the limit with $R\to\infty$ of the
integral over a contour as shown in figure \ref{tutaj} with $\tau=1+\delta$.
We can divide the contour into three parts: first starting at $-\infty$ and
ending at $c-i\tau$, the second from $c-i\tau$ to $R-i\tau$ ant the third part
from $R-i\tau$ to $R$ (see Figure \ref{tutaj}). As we have pointed out in the 
proof of Statement (4) the integral over the third part goes o $0$ as 
$R\to\infty$. Let $M_1$ be the value of the integral over the first part of
the contour. Then elementary computations show that there is a constant $M_2$
such that the integral in \refeq{modul} equals
\[
M_1+M_2|q^{-1}\y q^{is}|^{\tau}.
\]
Since $\tau=1+\delta$ there exist constants $C_1$ and $C_2$ such that
\[
\bigl|\an{\Phiq}{\y}{s-i\sigma}\bigr|\leq C_1+C_2
\bigl|q^{-1}\y q^{is}\bigr|^{1+\delta}
\]
for all $0\leq\sigma<1$ and all $s\in\RR$.
\qed

As an immediate consequence of Propositions \ref{Phiq} and \ref{uniq} we get
the following:
 
\begin{Cor}\label{Four1}
The functions\/ $\Fq$ and\/ $\Phiq$ are proportional:
\[
\Fq=-2\pi\rhbar\Fq(-q^2)^{-1}\Phiq.
\]
Moreover we have
\[
\Fq(\y)=-\frac{\Fq(-q^2)}{2\pi\rhbar}\Int_{\Yq}
\frac{\chi(-q^{-2}\y,\y')}
{\left(1-\overline{\y'}\right)\Fq(-q^2\y')}\,d\Mu(\y')
\]
in the sense of distribution theory (the correction of the integration contour
is understood as a part of the definition of the distribution under the sign
of the integral).
\end{Cor}

\begin{Cor}\label{wazne}
The distributional inverse Fourier transforms\/ $\widehat{\Fq}$ and\/
$\widehat{\overline{\Fq}}$ of\/ $\Fq$ and\/ $\overline{\Fq}$ satisfy
\[
\widehat{\Fq}(\y)=
\overline{\alpha(\y)}|\y|\chi(-1,\y)\widehat{\overline{\Fq}}(\y).
\]
\end{Cor}

\proof
By Corollary \ref{Four1}
\[
\widehat{\Fq}(\y)=-\frac{\Fq(-q^{2})}{2\pi\rhbar}
\frac{\chi(-q^{-2},\y)\y}{\left(1-\overline{\y}\right)\Fq(-q^2\y)}.
\]
Consequently
\[
\widehat{\overline{\Fq}}(\y)=\overline{\widehat{\Fq}(\y^{-1})}
=-\frac{\Fq(-q^2)^{-1}}{2\pi\rh}
\frac{\chi(-q^{-2},\y)(1-\y^{-1})}{\Fq(-q^{2}\y^{-1})\overline{\y}}.
\]
Therefore
\[
\frac{\widehat{\overline{\Fq}}(\y)}{\widehat{\Fq}(\y)}
=\Fq(-q^2)^{-2}\,\frac{\rhbar}{\rh}\,
\frac{1-\overline{\y}}{\y-1}
\Fq(-q^2\y)\Fq(-q^2\y^{-1})\,{\overline{\y}}^{-1}.
\]
Using \refeq{wzorek}, \refeq{FFal}, \refeq{AlfaChi}, the fact that $\chi$ is
a bicharacter and \refeq{bchar} we get
\[
\begin{array}{r@{\;=\;}l@{\smallskip}}
\displaystyle{
\frac{\widehat{\overline{\Fq}}(\y)}{\widehat{\Fq}(\y)}}
&\displaystyle{
-\Fq(-q^2)^{-2}\,\frac{\rhbar}{\rh}\,\Fq(-\y)\Fq(-q^2\y^{-1})\,
{\overline{\y}}^{-1}}\\
&\displaystyle{
-\Fq(-q^2)^{-2}\,\frac{\rhbar}{\rh}\,\Cq\alpha(-q^{-1}\y)\,
{\overline{\y}}^{-1}}\\
&\displaystyle{
-\Fq(-q^2)^{-2}\,\frac{\rhbar}{\rh}\,\Cq\alpha(-q^{-1})\alpha(\y)
\frac{\alpha(-q^{-1}\y)}{\alpha(-q^{-1})\alpha(\y)}\,{\overline{\y}}^{-1}}\\
&\displaystyle{
-\Fq(-q^2)^{-2}\,\frac{\rhbar}{\rh}\,\Cq\alpha(-q^{-1})\alpha(\y)
\chi(-q^{-1},\y)\,{\overline{\y}}^{-1}}\\
&\displaystyle{
-\Fq(-q^2)^{-2}\,\frac{\rhbar}{\rh}\,\Cq\alpha(-q^{-1})\alpha(\y)
\frac{1}{(\chi(-q,\y)\,{\overline{\y}}}}\\
&\displaystyle{
-\Fq(-q^2)^{-2}\,\frac{\rhbar}{\rh}\,\Cq\alpha(-q^{-1})\alpha(\y)
\frac{1}{\chi(q,\y){\overline{\y}}\chi(-1,\y)}}\\
&\displaystyle{
-\Fq(-q^2)^{-2}\,\frac{\rhbar}{\rh}\,\Cq\alpha(-q^{-1})\alpha(\y)
\frac{1}{\PHase{\y}{\overline{\y}}\chi(-1,\y)}}\\
&\displaystyle{
-\Fq(-q^2)^{-2}\,\frac{\rhbar}{\rh}\,\Cq\alpha(-q^{-1})
\frac{1}{\overline{\alpha(\y)}|\y|\chi(-1,\y)}}.
\end{array}
\]
Now inserting $\y=-1$ in \refeq{FFal} and using \refeq{wzorek} we obtain
\[
\Cq\alpha(-q^{-1})=\Fq(-1)\Fq(-q^2)=\frac{-\rh}{\rhbar}\Fq(-q^2)^2,
\]
and thus
$\widehat{\Fq}(\y)=\overline{\alpha(\y)}|\y|\chi(-1,\y)
\widehat{\overline{\Fq}}(\y)$.
\qed

\subsection{Other useful functions}

In this subsection we shall relate the special functions introduced above to
other functions which fit well the framework developed in \cite{qef}. This
step is needed to be able to freely use the theory of Zakrzewski relation
presented in that reference.

Consider the function $\PH\colon\Yq\to\TT$ given by
$\pH{q^nq^{it}}=e^{\frac{2\pi i n}{N}}$, where $N$ is the constant entering
the definition of $\rh$ (cf.~beginning of this section). It is easy to check
that this formula defines a function on $\Yq$ (the value does not depend on
the representation of $\y$ in the form $q^nq^{it}$). 

Elementary computations give the following formula
\begin{equation}\label{elem}
\y=\pH{\y}|\y|^{1+\frac{2\pi i\rre{\rh}}{N}}
\end{equation}  
for any $\y\in\Yq$.

In \cite{qef} for a parameter $-\pi<\hbar<\pi$ a subset $\Omega_\hbar^+$ was
defined in the following way
\[
\Omega_\hbar^+=\bigl\{z\in\CC\bez\{0\}:\:\hbar\,\mathrm{arg}\,z\geq0,\:
|\mathrm{arg}\,z|\in[0,|\hbar|]\bigr\}. 
\]

In what follows we shall use the identification of the region lying between
the logarithmic spirals $q^n\{q^{it}:\:t\in\RR\}$ and 
$q^{n+1}\{q^{it}:\:t\in\RR\}$
\[
\bigl\{q^nq^{i(t-i\tau)}:\:t\in\RR,\:\tau\in[0,1]\bigr\}
\]
with $\Omega_{\im{\rhinv}}^+$ given by
\[
q^nq^{i(t-i\tau)}\longleftrightarrow e^{i\tau\hbar}|q^nq^{it}|
\]
(notice that the conditions imposed on the form of $\rh$ imply that
$|\im{\rhinv}|<\pi$).

Another tool used in \cite{qef} is the space
\[
{\mathcal H}_\hbar^+=\left\{f\in C\left(\Omega_\hbar^+\right):\:
\begin{array}{c@{\smallskip}}
f\text{ is holomorphic in the interior of }\Omega_\hbar^+\text{ and}\\
\text{for any }\lambda>0\text{ the function }z\mapsto
e^{-\lambda\ell(z)^2}f(z)\\
\text{is bounded on }\Omega_\hbar^+
\end{array}
\right\},
\]
where $\ell(z)=\log{|z|}+i\,\mathrm{arg}\,z$. 

For $k\in\{0,\ldots,N-1\}$ and $r>0$ define 
\[
f_k(r)=\Fq\left(e^{\frac{2\pi ik}{N}}r^{1+\frac{2\pi i\rre{\rh}}{N}}\right)
\]

Using the information about the function $\Fq$ obtained in Theorem
\ref{glowne1} and Proposition \ref{wzrostFq} one can without difficulty prove
the following (cf.~\cite[Section 5.5]{phd} and \cite{azb}):

\begin{Lem}\label{male}
Let $k\in\{0,\ldots,N-1\}$. Then
\begin{itemize}
\item[{\rm (1)}] the function $f_k$ extends to a continuous function on
$\Omega_{\im{\rhinv}}^+$ which is holomorphic in the interior of that region,
moreover, denoting the extension by the same symbol, we have 
$f_k\in{\mathcal H}_{\im{\rhinv}}^+$ and for any $\tau\in]0,1[$ we have
$f_k^{-1}\in{\mathcal H}_{\tau\iim{\rhinv}}^+$;
\item[{\rm (2)}] for any $r>0$ we have
\[
f_k\left(e^{i\iim{\rhinv}}r\right)=(1+q^{-1}\y)\Fq(q^{-1}\y),
\]
where $\y=e^{\frac{2\pi ik}{N}}r^{1+\frac{2\pi i\rre{\rh}}{N}}$;
\item[{\rm (3)}] we have the following asymptotic behaviour of $f_k$:
\[
|f_k(z)|\leq\Theta_k(z)|q^{-1}z|^{\frac{\mathrm{arg}\,z}{\im{\rhinv}}},
\]
where $\lim\limits_{|z|\to\infty}\Theta_k(z)=1$.
\end{itemize}
\end{Lem}

\Section{Operator equalities}\label{opeq}

\subsection{The commutation relations}\label{komrel}

In this section we shall examine pairs $(A,B)$ of operators on a Hilbert space
satisfying the following conditions:
\begin{enumerate}
\item\label{c1} $A$ and $B$ are normal,
\item $\ker{A}=\ker{B}=\{0\}$,
\item $\spec{A},\:\spec{B}\subset\Ybar$,
\item\label{c4} for all $\y,\y'\in\Yq$
\begin{equation}\label{Weyl}
\chi(\y,A)\chi(B,\y')=\chi(\y,\y')\chi(B,\y')\chi(\y,A).
\end{equation}
\end{enumerate}
Formula \refeq{Weyl} is called the {\em Weyl relation}.
The set of pairs $(A,B)$ of operators on a Hilbert space $H$ fulfilling
conditions \ref{c1}--\ref{c4} will be denoted by $\Dd_H$. 

For any infinite dimensional Hilbert space $H$ the set $\Dd_H$ is non empty.
Moreover any pair $(A,B)\in\Dd_H$ is unitarily equivalent to a direct sum of
so called {\em Schr\"odinger pairs}. The Schr\"odinger pair 
$(A_{\text{\tiny S}},B_{\text{\tiny S}})$ acts irreducibly on $L^2(\Yq)$ in
the following way:
\[
\begin{array}{r@{\;=\;}l@{\smallskip}}
\bigl(A_{\text{\tiny S}}f\bigr)(\y)&\an{f}{q\y}{-i},\\
\bigl(B_{\text{\tiny S}}f\bigr)(\y)&\y f(\y)
\end{array}
\]
(cf.~\cite{azb}, \cite{phd}).

The correspondence $H\mapsto\Dd_H$ satisfies the following conditions:
\begin{itemize}
\item If $H$ and $K$ are Hilbert spaces, $U\colon H\to K$ a unitary map and
$(A,B)\in\Dd_H$ then $(UAU^*,UBU^*)\in\Dd_K$,
\item if $H=H_1\oplus H_2$ and operators $A$ and $B$ on $H$ decompose as
$A=A_1\oplus A_2$ and $B=B_1\oplus B_2$ respectively then $(A,B)\in\Dd_H$ if
and only if $(A_k,B_k)\in\Dd_{H_k}$ for $k=1,2$.
\end{itemize}
Such a structure is called an {\em operator domain}, a notion closely related
to {\em compact}\/ and {\em measurable domains}\/ and {\em
$\mathrm{W}^*$-categories}\/ (cf.~\cite{ncgnt,dualth,gen,azb,phd}).
We shall use the symbol $\Dd$ also to denote this operator domain.

Inserting in the Weyl relation $(q,q^{it})$, $(q^{it},q)$, $(q,q)$, and
$(q^{it},q^{it'})$ for $(\y,\y')$ and
performing analytic continuation with respect to $t$ in the first two
cases we obtain (cf.~\refeq{bchar}):
\begin{equation}\label{finalrel}
\begin{array}{r@{\;=\;}l}
\Phase{A}|B|&|q||B|\Phase{A},\\
|A|\Phase{B}&|q|\Phase{B}|A|,\\
\Phase{A}\Phase{B}&e^{i\iim{\rhinv}}\Phase{B}\Phase{A},\\
|A|^{it}|B|^{it'}&e^{-itt'\iim{\rhinv}}|B|^{it'}|A|^{it},
\end{array}
\end{equation}
The last equation of \refeq{finalrel} means that $|B|$ and $|A|$ satisfy the
{\em Zakrzewski relation}\/ (\cite[Definition 2.1]{qef}) with 
$\hbar=\im{\rhinv}$. We
shall represent this graphically as $\Zrel{|B|}{|A|}{\hbar}$. 

Let us recall the definition of a {\em core for a family of operators}\/
introduced in \cite{gl}. Let $H$ be a Hilbert space and let $\mathscr{T}$ be a
family of closable operators on $H$. A linear subset $D_0$ is a {\em core for
$\mathscr{T}$}\/ if 
\begin{enumerate}
\item $D_0\subset\D{T}$ for all $T\in\mathscr{T}$,
\item for any $x\in H$ there exists a sequence $(x_n)_{n\in\NN}$ of elements
of $D_0$ converging to $x$ such that
\[
\left(\begin{array}{c}T\in\mathscr{T},\\ x\in\D{\overline{T}}\end{array}\right)
\Longrightarrow\Bigl(\:Tx_n\xrightarrow[n\to\infty]{}\overline{T}x\:\Bigl).
\]
\end{enumerate}

Throughout the paper we shall adopt the following convention. For any pair of 
linear operators $(X,Y)$ acting on the same Hilbert space we shall denote by
$X\comp Y$ their composition. If $X\comp Y$ happens to be closable then $XY$
will denote the closure of $X\comp Y$.

Let us recall \cite[Theorem 2.7]{gl}.

\begin{Thm}\label{corethm}\sloppy
Let $H$ be a Hilbert space and let $T_1,\ldots,T_p$ be normal operators on
$H$. 
Assume that for each pair of
indices $k,l\in\{1,\ldots,p\}$ there exist scalars $\mu(k,l)$
and $\lambda(k,l)$ such that
\[
\begin{array}{r@{\;=\;}l@{\smallskip}}
\Phase{T_k}|T_l|&\mu(k,l)|T_l|\Phase{T_k},\\
\Phase{T_k}\Phase{T_l}&\lambda(k,l)\Phase{T_l}\Phase{T_k}.
\end{array}
\]
Moreover
assume that there exists a real number $\hbar$ with $|\hbar|<\pi$ such that
for any $k,l\in\{1,\ldots,p\}$ one of the following conditions holds:
\begin{enumerate}
\item $|T_k|$ strongly commutes with $|T_l|$,
\item $\Zrel{|T_k|}{|T_l|}{\hbar}$,
\item $\Zrel{|T_l|}{|T_k|}{\hbar}$.
\end{enumerate}
Let $\mathscr{T}$ be the family of all compositions of the form 
$T_{i_1}^\sharp \comp T_{i_2}^\sharp\comp\cdots\comp T_{i_n}^\sharp$ where
$i_1,\ldots,i_n\in\{1,\ldots,p\}$ and $T^\sharp$ denotes either $T$ or $T^*$. 
Then
all operators in $\mathscr{T}$ are densely defined and closable and there
exists a core for $\mathscr{T}$.
\end{Thm}

The original formulation of this theorem in \cite{gl} places stronger
conditions on the operators $T_1,\ldots,T_p$. Namely the constant $\hbar$ is
fixed as $\frac{2\pi}{N}$ with an even natural number $N$ greater or equal to
$6$ and the numbers $\lambda(k,l),\:\mu(k,l)$  are equal to $1$ for all
$k,l\in\{1,\ldots,p\}$. It is, however, clear from the proof given in
\cite{gl} that these restrictions are not essential. 

Given a pair $(A,B)\in\Dd_H$ for some Hilbert space $H$ we can apply Theorem 
\ref{corethm} with $p=2$ and $T_1=A$, $T_2=B$. The family of all finite
compositions of elements of the set $\{A,A^*,B,B^*\}$ will be denoted by
$\mathscr{T}$ and we shall use the symbol $D_0$ for the core for $\mathscr{T}$. 

\subsection{Products}

\begin{Thm}\label{iloczyny}
Let $H$ be a Hilbert space and let $(A,B)\in\Dd_H$. Then
\begin{itemize}
\item[{\rm (1)}] The operators $A\comp B$, $B\comp A$, 
$A\comp B^*$ and $B^*\comp A$
are closable and we have
\begin{equation}\label{produkty}
AB=q^2BA,\qquad AB^*=B^*A;
\end{equation}
\item[{\rm (2)}] for any $\y\in\Yq$
\begin{equation}\label{alphachi}
\alpha(\y)\chi(B,\y)\chi(\y,A)=\alpha(B)^*\chi(\y,A)\alpha(B)
=\alpha(A)\chi(B,\y)\alpha(A)^*;
\end{equation}
\item[{\rm (3)}] we have
\begin{equation}\label{unit_rown}
qBA=\alpha(B)^*A\alpha(B)=\alpha(A)B\alpha(A)^*.
\end{equation}
\end{itemize}
\end{Thm}

\proof
\sloppy{\sc Ad (1).}
The closability of all finite compositions of elements of the set
$\{A,A^*,B,B^*\}$ was established in Theorem \ref{corethm}. Formula
\refeq{produkty} follows from the fact that $A(By)=q^2B(Ay)$ and
$A(B^*y)=B^*(Ay)$ for any $y\in D_0$.

\sloppy{\sc Ad (2).}
In order to make our exposition shorter we shall only prove the equality
$\alpha(\y)\chi(B,\y)\chi(\y,A)=\alpha(A)\chi(B,\y)\alpha(A)^*$. The other one
can be proved in an analogous fashion.

From the Weyl relation and the fact that $\chi$ is a bicharacter we infer that
\begin{equation}\label{przeksz1}
\chi(B,\y)^*\chi(\y',A)\chi(B,\y)=\chi(\y',\y A).
\end{equation}
Inserting $\y'=q$ and $\y'=q^{it}$ in \refeq{przeksz1} we obtain
\[
\begin{array}{r@{\;=\;}l@{\smallskip}}
\chi(B,\y)^*\Phase{A}\chi(B,\y)&\Phase{\y}\Phase{A},\\
\chi(B,\y)^*|A|\chi(B,\y)&|\y||A|,
\end{array}
\]
where in the second equality we performed analytic continuation to $t=-i$.
Multiplying left and right hand sides of these equations results in
\begin{equation}\label{tualpha}
\chi(B,\y)^*A\chi(B,\y)=\y A
\end{equation}
for all $\y\in\Yq$. Now let us apply function $\alpha$ to both sides of
\refeq{tualpha} and multiply both sides of the resulting equation by
$\alpha(A)^*$ from the right. We obtain
\[
\alpha(A)\chi(B,\y)\alpha(A)^*=\chi(B,\y)\alpha(\y A)\alpha(A)^*.
\]
Remembering that $\alpha(\y)\overline{\alpha(\y)}=1$ for all $\y\in\Yq$ we can
rewrite this as
\[
\begin{array}{r@{\;=\;}l@{\smallskip}}
\alpha(A)\chi(B,\y)\alpha(A)^*
&\alpha(\y)\chi(B,\y)\bigl(\alpha(\y A)\overline{\alpha(\y)}\alpha(A)^*\bigr)\\
&\alpha(\y)\chi(B,\y)
\bigl(\alpha(\y A)\left[\alpha(\y)\alpha(A)\right]^{-1}\bigr).
\end{array}
\]
Our assertion follows now from \refeq{AlfaChi}.

{\sc Ad (3).}
As in the proof of (2) w put $\y=q$ and $\y=q^{it}$ in \refeq{alphachi} 
and perform analytic continuation in the latter case:
\[
\begin{array}{r@{\;=\;}l@{\smallskip}}
e^{\frac{i}{2}\iim{\rhinv}}\Phase{B}\Phase{A}
&\alpha(B)^*\Phase{A}\alpha(B)=\alpha(A)\Phase{B}\alpha(A)^*,\\
e^{\frac{i}{2}\iim{\rhinv}}|B||A|&\alpha(B)^*|A|\alpha(B)
=\alpha(A)|B|\alpha(A)^*.
\end{array}
\]
Now we can multiply the left and right hand sides of the above equations and
use relations \refeq{finalrel} to obtain \refeq{unit_rown}.
\qed

\begin{Cor}\label{corAB}
Let $H$ be a Hilbert space and let $(A,B)\in\Dd_H$. Then 
\begin{itemize}
\item[{\rm (1)}] for any $\y_1,\y_2\in\Yq$ we have $(\y_1A,\y_2B)\in\Dd_H$;
\item[{\rm (2)}] $(B,A^{-1})\in\Dd_H$
\item[{\rm (3)}] $(AB,B),(A,BA)\in\Dd_H$;
\item[{\rm (4)}] for any $\y\in\Yq$ we have 
$\chi(q^{-1}AB,\y)=\alpha(\y)\chi(B,\y)\chi(\y,A)$.
\end{itemize}
\end{Cor}

\proof
{\sc Ad (1).}
The assertion follows by multiplying both sides of the Weyl relation
\ref{Weyl} by $\chi(\y_2,\y_1)$ and using the fact that $\chi$ is a
bicharacter.

{\sc Ad (2).} Applying hermitian conjugation to both sides of the Weyl 
relation and using the fact that 
\[
\begin{array}{r@{\;=\;}l@{\smallskip}}
\chi(\y,A)^*&\chi(\y,A^{-1}),\\
\chi(B,\y')^*&\chi(B,{\y'}^{-1})
\end{array}
\]
we obtain
\[
\chi(B,{\y'}^{-1})\chi(\y,A^{-1})=\chi({\y'}^{-1},\y)
\chi(\y,A^{-1})\chi(B,{\y'}^{-1})
\]
which by symmetry of $\chi$ means that $(B,A^{-1})\in\Dd_H$.

{\sc Ad (3).} This follows by conjugating the pair $(A,B)$ with the unitary
operators $\alpha(B)$ and $\alpha(A)^*$ and using Statement (3) of 
Theorem \ref{iloczyny} and Statement (1) above.

{\sc Ad (4).} By Statement (1) of Theorem \ref{iloczyny} we have
\[
\alpha(\y)\chi(B,\y)\chi(\y,A)=\alpha(A)\chi(B,\y)\alpha(A^*)
=\chi(\alpha(A)B\alpha(A)^*,\y).
\]
Now Statement (2) of the same theorem says that this is equal to
\[
\chi(qBA,\y)=\chi(q^{-1}AB,\y).
\]
\qed

\subsection{Sums}\label{SumSubs}

Let $H$ be a Hilbert space and let $(A,B)\in\Dd_H$. Put $S=\pH{A}|A|^{\rh}$ and
$R=\pH{B}|B|^{\rh}$. Then the pair $(R,S)$ satisfies the commutation relations
and spectral conditions considered in \cite{azb}. Therefore we can use 
Proposition 2.3 of that reference which yields

\begin{Prop}\label{gRupa}
Let $H$ be a Hilbert space and let $(A,B)\in\Dd_H$. Then there exists a one
parameter group $({\mathcal R}_t)_{t\in\RR_+}$ of unitary operators acting on
$H$ such that 
\[
\begin{array}{r@{\;=\;}l@{\qquad}r@{\;=\;}l@{\smallskip}}
\pH{B}\mathcal{R}_t&\mathcal{R}_t\pH{B},
&\pH{A}\mathcal{R}_t&\mathcal{R}_t\pH{A},\\
\mathcal{R}_t|B|&|B|^t\mathcal{R}_t,
&|A|\mathcal{R}_t&\mathcal{R}_t|A|^t.
\end{array}
\]  
for all $t\in\RR_+$.
\end{Prop}

In the next theorem we shall consider the operator 
$A+B\comp A$, where $(A,B)\in\Dd_H$ for some Hilbert space $H$. It turns out
to be closable and its closure coincides with the closure of $A+BA$ (which, of
course, is closable as well). The closability is more or less straightforward:
$(A+B\comp A)^*\subset A^*+A^*\comp B^*$. The core $D_0$ described after
Theorem \ref{corethm} is contained in the domain of the latter operator which
is thus densely defined. The same reasoning shows that $A+BA$ is closable.

Let $x\in\D{A}\cap\D{BA}$ and let $(x_n)_{n\in\NN}$ be a sequence of elements
of $D_0$ such that
\[
\begin{array}{r@{\;\xrightarrow[n\to\infty]{}\;}l@{\smallskip}}
x_n&x,\\ Ax_n&Ax,\\ BAx_n& BAx.
\end{array}
\]
It follows that $\bigl((A+B\comp A)x_n\bigr)_{n\in\NN}$ is convergent and
consequently $x$ belongs to the domain of the closure of $A+B\comp A$.
Moreover
\[
\left(\overline{A+B\comp A}\right)x=\lim_{n\to\infty}(A+B\comp A)x_n
=\lim_{n\to\infty}(Ax_n+BAx_n)=Ax+BAx=\left(\overline{A+BA}\right)x.
\]
As $\D{A}\cap\D{BA}$ is a core for $\overline{A+BA}$ we obtain
$\overline{A+BA}\subset\overline{A+B\comp A}$.
The converse inclusion is trivial and our assertion follows. From now on the
closure of a sum of operators will be denoted by the symbol ``$\dplus$''. For
example
\[
\overline{A+B\comp A}=\overline{A+BA}=A\dplus BA.
\]

The proof of the theorem below is almost identical to proofs of analogous
theorems in \cite{OpEq} and \cite{azb}. Since the details are somewhat
different and our notation is not fully compatible with the one used in these
references, we decided to include this proof for the reader's convenience. 

\begin{Thm}\label{sumthm}
Let $H$ be a Hilbert space and let $(A,B)\in\Dd_H$. Then
\begin{itemize}
\item[{\rm (1)}] the operator $A+BA$ is densely defined and closable and its 
closure
\begin{equation}\label{SplRS}
A\dplus BA=\Fq(B)^*A\Fq(B),
\end{equation}
in particular $A\dplus BA$ is normal and\/ $\Spec{A\dplus BA}\subset\Ybar$;
\item[{\rm (2)}] the operator $A+B$ is densely defined and closable and its 
closure
\begin{equation}\label{suma1}
A\dplus B=\Fq(BA^{-1})^*A\Fq(BA^{-1}),
\end{equation} 
and
\begin{equation}\label{suma2}
A\dplus B=\Fq(B^{-1}A)B\Fq(B^{-1}A)^*;
\end{equation}
it follows that $A\dplus B$ is normal and\/ $\Spec{A\dplus B}\subset\Ybar$;
\item[{\rm (3)}] the function $\Fq$ has the {\em exponential property:}
\begin{equation}\label{expon}
\Fq(A\dplus B)=\Fq(B)\Fq(A).
\end{equation}
\end{itemize}
\end{Thm}

\proof
{\sc Ad (1).} We already know that $A+BA$ is closable. First we shall show
that $A\dplus BA\subset\Fq(B)^*A\Fq(B)$. 

It is easy to see that $\pH{B}^N=I$, so that the spectrum of $\pH{B}$ is
contained in the set of $N^{\text{th}}$ roots of unity. Let 
\begin{equation}\label{rozkladH}
H=\bigoplus_{k=0}^{N-1}
\end{equation}
be the decomposition into eigenspaces of $\pH{B}$. It is also easy to check
(cf.~remarks preceding Proposition \ref{gRupa}) that $|A|$ commutes with
$\pH{B}$ and thus preserves the decomposition \refeq{rozkladH}. We can
therefore consider each summand of this decomposition separately. 

We know that $\Zrel{|B|}{|A|}{\hbar}$ with $\hbar=\im{\rhinv}$. By
\cite[Theorem 3.1(2)]{qef} on each subspace $H_k$ we have
\[
f_k\left(e^{i\iim{\rhinv}}|B|\right)\subset|A|f_k(|B|),
\]
which by Lemma \ref{male}(2) means that
\[
\left(1+q^{-1}e^{\frac{2\pi i}{N}}|B|^{1+\frac{2\pi i\rre{\rh}}{N}}\right)
\Fq\left(q^{-1}e^{\frac{2\pi i}{N}}|B|^{1+\frac{2\pi i\rre{\rh}}{N}}\right)|A|
\subset|A|\Fq\left( q^{-1}e^{\frac{2\pi i}{N}}
|B|^{1+\frac{2\pi i\rre{\rh}}{N}}\right).
\]
Taking a direct sum over $k$ and keeping in mind \refeq{elem} we get 
\[
(1+q^{-1}B)\Fq(q^{-1}B)|A|\subset|A|\Fq(B)
\]
which can be rewritten as
\begin{equation}\label{zaw1}
|A|+q^{-1}B\comp|A|\subset\Fq(q^{-1}B)^*|A|\Fq(B).
\end{equation}
Finally multiplying both sides of \refeq{zaw1} from the left by $\phase{A}$, 
using \refeq{finalrel} and taking closure of both sides we obtain
\[
A\dplus BA\subset\Fq(B)^*A\Fq(B).
\]

In order to see the converse inclusion we shall prove that $\D{A+B\comp A}$ is
a core for $\Fq(B)^*A\Fq(B)$. Notice that apart from the relation 
$\Zrel{|B|}{|A|}{\im{\rhinv}}$ (cf.~\refeq{finalrel} and following remarks), 
for any $\tau\in]0,1[$ we have 
$\Zrel{|B|}{|A|^{\tau}}{\tau\iim{\rhinv}}$. Therefore by Statement (1) of 
Lemma \ref{male} and \cite[Theorem 3.1(3)]{azb} we have
\begin{equation}\label{rownosc}
f_k\left(e^{i\tau\iim{\rhinv}}|B|\right)|A|^\tau=|A|^\tau f_k(|B|)
\end{equation}
on $H_k$. Let $x\in\D{\Fq(B)^*A\Fq(B)}$ and let $x_k$ be the projection of $x$
onto $H_k$. Since for any $\tau\in]0,1[$ we have $\Fq(B)x\in\D{|A|^\tau}$ and
$\pH{B}$ commutes with $|A|$, the vectors $\Fq(B)x_k\in\D{|A|^\tau}$. Moreover
on $H_k$ we have $\Fq(B)=f_k(|B|)$. Therefore
\begin{equation}\label{Drozklad}
f_k(|B|)x_k\in\D{|A|^\tau}.
\end{equation}
Comparing \refeq{rownosc} and \refeq{Drozklad} yields
$x_k\in\D{f_k\left(e^{i\tau\iim{\rhinv}}|B|\right)|A|^\tau}$,
which by Proposition \ref{gRupa} gives
\begin{equation}\label{impl}
\mathcal{R}_\tau x_k\in
\D{f_k\left(e^{i\tau\iim{\rhinv}}|B|^{\tau^{-1}}\right)|A|}.
\end{equation}

By Statement (1) of Lemma \ref{male} the function 
$\RR_+\ni r\mapsto
\left|f_k\left(e^{i\tau\iim{\rhinv}}r^{\tau^{-1}}\right)\right|$ is separated
from $0$ and by Statement (3) of that lemma, behaves asymptotically like the 
function $r\mapsto r$. Therefore \refeq{impl} implies
that $|A|\mathcal{R}_\tau x_k\in\D{B}$ and consequently
$A\mathcal{R}_\tau x_k\in\D{B}$ ($\phase{A}$ preserves the domain of $B$).
Taking the direct sum over $k$ we see that for $\tau\in]0,1[$ 
\[
\mathcal{R}_\tau x\in\D{A},\qquad A\mathcal{R}_\tau x\in\D{B}.
\]
It follows that 
\begin{equation}\label{Dsumy}
\mathcal{R}_\tau x\in\D{A+B\comp A}.
\end{equation}
Moreover
\[
\begin{array}{r@{\;}c@{\;}l@{\smallskip}}
\Fq(B)^*A\Fq(B)\mathcal{R}_\tau x&=&
\Fq\left(\pH{B}|B|^{1+\frac{2\pi i\rre{\rh}}{N}}\right)^*A\Fq\left(
\pH{B}|B|^{1+\frac{2\pi i\rre{\rh}}{N}}\right)\mathcal{R}_\tau x\\
&=&\mathcal{R}_\tau\Fq\left(\pH{B}|B|^{\tau^{-1}\left(1
+\frac{2\pi i\rre{\rh}}{N}\right)}\right)^*
\pH{A}|A|^{\tau\left(1+\frac{2\pi i\rre{\rh}}{N}\right)}\\
&&\times\Fq\left(\pH{B}|B|^{\tau^{-1}\left(1
+\frac{2\pi i\rre{\rh}}{N}\right)}\right)x.
\end{array}
\]
As the group $(\mathcal{R}_t)_{t\in\RR_+}$ is strongly continuous we conclude
that
\[
\begin{array}{r@{\;\xrightarrow[\tau\nearrow 1]{}\;}l@{\smallskip}}
\mathcal{R}_\tau x&x,\\
\Fq(B)^*A\Fq(B)\mathcal{R}_\tau x&\Fq(B)^*A\Fq(B)x.
\end{array}
\]
Now \refeq{Dsumy} implies that $\D{A+B\comp A}$ is a core for
$\Fq(B)^*A\Fq(B)$.

\begin{Rem}\label{Fakcik} 
Before continuing the proof Let us notice a fact which we
established in the proof of Statement (1):
\medskip

\noindent{\em
Let $H$ be a Hilbert space and let $(A_0,B_0)\in\Dd_H$. Then there is a 
one parameter group $(\mathcal{R}_t)_{t\in\RR_+}$ of unitary operators acting
on $H$ such that for all $t\in\RR_+$ 
\begin{equation}\label{Regu}
\begin{array}{r@{\;=\;}l@{\smallskip}}
\mathcal{R}_t|B_0|\mathcal{R}_t^*&|B_0|^t,\\
\mathcal{R}_t^*|A_0|\mathcal{R}_t&|A_0|^t
\end{array}
\end{equation}
and for any $x\in\D{A_0\dplus B_0A_0}$ and any $t>0$
the vector $\mathcal{R}_tx$ belongs to $\D{A_0}\cap\D{B_0\comp A_0}$ and the net
$(\mathcal{R}_tx)_{t]0,1[}$ converges to $x$ in the graph topology of
$A_0\dplus B_0A_0$ as $t\nearrow1$.}
\end{Rem}\medskip
     
{\sc Ad (2).} By repeated use of Corollary \ref{corAB} we conclude that
$(A,BA^{-1})\in\Dd_H$. Formula \ref{SplRS} applied to this pair yields
\refeq{suma1}.

To prove \refeq{suma2} let us observe that by Theorem \ref{FFALT} we have
\[
\Fq(\y)=\Cq\alpha(q^{-1}\y)\overline{\Fq(q^2\y^{-1})}
\]
for all $\y\in\Yq$. Therefore
\begin{equation}\label{subst}
\Fq(BA^{-1})=\Cq\alpha(q^{-1}BA^{-1})\Fq(q^{2}AB^{-1})^*.
\end{equation}
Moreover the pair $(A,q^{-1}BA^{-1})$ belongs to $\Dd_H$, so by Statement (3)
of Theorem \ref{iloczyny}
\begin{equation}\label{transfal}
\alpha(q^{-1}BA^{-1})^*A\alpha(q^{-1}BA^{-1})=B.
\end{equation}
Inserting \refeq{subst} into \refeq{suma1} and using \refeq{transfal}
we obtain \refeq{suma2}.

{\sc Ad (3).}
Let $T=B^{-1}A$. By Corollary \ref{corAB} the pairs $(T,B),(T,A)\in\Dd_H$.
Then applying the transformation $m\mapsto\Fq(B^{-1}A)m\Fq(B^{-1}A)^*$ to the
pair $(T,B)$ we get $(T,A\dplus B)\in\Dd_H$ by \refeq{suma2}. Inserting this
last pair in place of $(A,B)$ in \refeq{suma1} we obtain
\begin{equation}\label{prawa1}
\Fq(A\dplus B)^*T\Fq(A\dplus B)=T\dplus(A\dplus B)T,
\end{equation}
while since $(T,B),(T,A)\in\Dd_H$ and $A$ strongly commutes with $BT$, we have
\begin{equation}\label{prawa2}
\begin{array}{r@{\;=\;}l@{\smallskip}}
\Fq(A)^*\Fq(B)^*T\Fq(B)\Fq(A)&\Fq(A)^*\left(T\dplus BT\right)\Fq(A)\\
&\Fq(A)^*T\Fq(A)\dplus BT=(T\dplus AT)\dplus BT.
\end{array}
\end{equation}

Notice that $T$ strongly commutes with $BA^{-1}=q^2T^{-1}$ and that
$(AT,BA^{-1})\in\Dd_H$ (as $(T,A)\in\Dd_H$ implies
$(T,AT)\in\Dd_H$ and consequently $(AT,BA^{-1})=(AT,T^{-1})\in\Dd_H$).
Therefore by \refeq{suma1} and \refeq{SplRS}
\begin{equation}\label{roz}
\begin{array}{r@{\;=\;}l@{\smallskip}}
(A\dplus B)T&\Fq(BA^{-1})^*A\Fq(BA^{-1})T\\
&\Fq(BA^{-1})^*AT\Fq(BA^{-1})\\
&\overline{AT+BA^{-1}\comp AT}=AT\dplus BT.
\end{array}
\end{equation}

In view of Equation \refeq{roz} the right hand side of \refeq{prawa1} is equal
to $T\dplus(AT\dplus BT)$. Moreover the operators on the right hand sides of
\refeq{prawa1} and \refeq{prawa2} coincide on the subspace 
$D(T)\cap D(AT)\cap D(BT)$.
 
We shall need the following

\begin{Lem}\label{Core1}
$D(T)\cap D(AT)\cap D(BT)$ is a core for $T\dplus(A\dplus B)T$.
\end{Lem}

From this Lemma we immediately see that 
$T\dplus(A\dplus B)T\subset(T\dplus AT)\dplus BT$ and since both these
operators are normal, we have $T\dplus(A\dplus B)T=(T\dplus AT)\dplus BT$.
It follows now from formulae \refeq{prawa1} and \refeq{prawa2} that the
operator $\Fq(B)\Fq(A)\Fq(A\dplus B)^*$ commutes with $T$ and thus with
$|T|^{it}$ for all $t\in\RR$:
\[
\Fq(B)\Fq(A)\Fq(A\dplus B)^*=|T|^{it}\Fq(B)\Fq(A)\Fq(A\dplus B)^*|T|^{-it}.
\]
On the other hand, since $(T,B)$, $(T,A)$ and $(T,(A\dplus B))$ belong to
$\Dd_H$, we have (cf.~\refeq{finalrel})
\[
\begin{array}{r@{\;=\;}l@{\smallskip}}
|T|^{it}B|T|^{-it}&q^{it}B,\\
|T|^{it}A|T|^{-it}&q^{it}A,\\
|T|^{it}(A\dplus B)|T|^{-it}&q^{it}(A\dplus B).
\end{array}
\]
Consequently
\begin{equation}\label{final}
\Fq(B)\Fq(A)\Fq(A\dplus B)^*=\Fq(q^{it}B)\Fq(q^{it}A)\Fq(q^{it}(A\dplus B))^*
\end{equation}
for all $t\in\RR$. The right hand side converges strongly to $I$ as
$t\to-\infty$ while the left hand side is independent of $t$. Therefore 
$\Fq(B)\Fq(A)\Fq(A\dplus B)^*=I$ and \refeq{expon} follows.
\qed

The formula \refeq{expon} justifies the name {\em quantum exponential
function}\/ used in Subsection \ref{qefun}.

\noindent{\sc Proof of Lemma \ref{Core1}.}
We shall first prove that $\D{T^2}\cap\D{(A\dplus B)T}$ is a core for
$T\dplus(A\dplus B)T$. 

Since $(T,A\dplus B)\in\Dd_H$ we have $\Zrel{|T|}{|A\dplus B|}{\hbar}$ with
$\hbar=-\im{\rhinv}$. For $\tau>0$ define
$f_\tau\colon\overline{\Omega_{-\im{\rhinv}}^+}\to\CC$ by
\[
f_\tau(z)=\left\{\begin{array}{c@{\quad\text{ for }}l}
e^{-\tau(\log{z})^2}&z\neq0,\\0&z=0.
\end{array}\right.
\]
Then $f_\tau$ is a bounded and continuous functions on
$\overline{\Omega_{-\im{\rhinv}}^+}$ which is holomorphic in the interior of
this set. Moreover $f_\tau$ converges almost uniformly on 
$\Omega_{-\im{\rhinv}}^+$ to the constant function equal to $1$, as
$\tau\searrow0$. By \cite[Theorem 3.1(2)]{qef} we have
\[
f_\tau\left(e^{-i\iim{\rhinv}}|T|\right)
|A\dplus B|\subset|A\dplus B|f_\tau(|T|),
\]
and it follows that
\[
f_\tau\left(e^{-i\iim{\rhinv}}|T|\right)
|A\dplus B|T\subset|A\dplus B|Tf_\tau(|T|).
\]
Multiplying both sides of this equation from the left by $\PHase{A\dplus B}$
and using relations \refeq{finalrel} (with $(T,A\dplus B)$ in place of
$(A,B)$) we obtain
\[
f_\tau(q^{-1}|T|)(A\dplus B)T\subset(A\dplus B)Tf_\tau(|T|).
\]
In particular for any $x\in\D{(A\dplus B)T}$ we have
$f_\tau(|T|)x\in\D{(A\dplus B)T}$. Also, as the function 
$z\mapsto|z|^2f_\tau(|z|)$ is bounded, the vector $x\in\D{T^2}$. By the almost
uniform convergence of $(f_\tau)_{\tau>0}$ the net
$f_\tau(|T|)x\xrightarrow[\tau\searrow0]{}x$ in the graph topology of 
$T\dplus(A\dplus B)T$. Indeed:
\[
\begin{array}{r@{\;=\;}l@{\smallskip}}
\bigl(T\dplus(A\dplus B)T\bigr)f_\tau(|T|)x&
Tf_\tau(|T|)x+(A\dplus B)Tf_\tau(|T|)x\\
&f_\tau(|T|)Tx+f_\tau(q^{-1}|T|)(A\dplus B)Tx
\end{array}
\]
and the right hand side converges to 
$Tx+(A\dplus B)Tx=\bigl(T\dplus(A\dplus B)T\bigr)x$. We have therefore proved
that $\D{T^2}\cap\D{(A\dplus B)T}$ is a core for
$T\dplus(A\dplus B)T$.

Now we shall use the fact stated in Remark \ref{Fakcik} with
$(A_0,B_0)=(AT,BA^{-1})$. Let 
$x\in\D{AT\dplus BT}=\D{\overline{AT+BA^{-1}\comp AT}}$. We have 
$\mathcal{R}_tx\in\D{AT}\cap\D{BT}$ and $\mathcal{R}_tx$ converges to $x$ in
the graph topology of $AT\dplus BT$, as $t\nearrow1$.

Since $BA^{-1}=q^2T^{-1}$, by \refeq{Regu} we see that
\begin{equation}\label{regu2}
\mathcal{R}_t^*|T|\mathcal{R}_t=|T|^{t^{-1}}
\end{equation}
for all $t\in\RR_+$. Assume now that $\frac{1}{2}\leq t<1$ and that
$x\in\D{T^2}$. We have $x\in\D{|T|^{t^{-1}}}$, because $t^{-1}\leq2$. 
Therefore by \refeq{regu2} the vectors $\mathcal{R}_tx\in\D{T}$. 
Now rewriting \refeq{regu2} as
\[
|T|\mathcal{R}_tx=\mathcal{R}_t|T|^{t^{-1}}x,
\]
shows that $\mathcal{R}_tx\xrightarrow[t\nearrow1]{}x$ in the graph topology 
of $T$. 

This way we showed that if $x\in\D{T^2}\cap\D{(A\dplus B)T}$ then 
$\mathcal{R}_tx\in\D{T}\cap\D{AT}\cap\D{BT}$ provided that
$\frac{1}{2}\leq t<1$ and moreover 
$\mathcal{R}_tx\in\D{T}\cap\D{AT}\cap\D{BT}$ and
$\mathcal{R}_tx\xrightarrow[t\nearrow1]{}x$ in the graph topologies of
$T$ and $AT\dplus BT$. It follows that $\mathcal{R}_tx$ converges to $x$ 
in the graph topology of $T\dplus(AT\dplus BT)=T\dplus(A\dplus B)T$. 

Now as we have established that $\D{T^2}\cap\D{(A\dplus B)T}$ is a core for 
$T\dplus(A\dplus B)T$, it follows that so is $\D{T}\cap\D{AT}\cap\D{BT}$.
\qed
     
\begin{Cor}
Let $H$ be a Hilbert space and let $(A,B)\in\Dd_H$. Then
\begin{equation}\label{F5}
\Fq(BA)=\Fq(B)^*\Fq(A)\Fq(B)\Fq(A)^*.
\end{equation}
\end{Cor}

\proof
Applying $\Fq$ to both sides of \refeq{SplRS} we obtain
\[
\Fq(A\dplus BA)=\Fq(B)^*\Fq(A)\Fq(B).
\]
On the other hand since $(A,BA)\in\Dd_H$, by \refeq{expon} we have
\[
\Fq(A\dplus BA)=\Fq(BA)\Fq(A).
\] 
Comparing these formulae yields \refeq{F5}.
\qed

\subsection{Necessity of the spectral condition}

In this subsection we would like to present another aspect of analysis of
the commutation relations we have considered so far. This aspect is not
relevant to the construction of new quantum ``$az+b$'' groups, so we have
decided to present it without proof. We refer to \cite[Section 6.5]{phd} for
details (cf.~also \cite[Section 2]{OpEq}).

One can consider pairs $(A,B)$ of normal 
operators on a Hilbert space $H$ satisfying 
more general commutation relations than those implied by the definition of
the operator domain $\Dd$ (we consider spectral conditions as part of the
commutation relations). Namely one can ask only that $(A,B)$ satisfy the
relations \refeq{finalrel}. Then it is easy to show that the spectra of $A$
and $B$ are contained in closures of unions of orbits of the group $\Yq$ in
$\CC$. If we assume that $(A,B)$ is irreducible (the only projections
commuting with $\PHase{A}$, $\PHase{B}$, $|A|^{it}$ and $|B|^{it}$ for all
$t\in\RR$ are $0$ and $I$) then the spectra are precisely equal to closures of
single orbits. Then multiplying both operators by a non zero scalar we can
suppose that $\spec{A}=\Ybar$. The spectrum of $B$ will coincide with
$\lambda\Ybar$ for some non zero $\lambda\in\CC$. 

It is clear that in the above situation $(A,B)\in\Dd_H$ if and only if
$\lambda\in\Yq$ or equivalently $\spec{B}=\Ybar$. It turns out that this
condition is equivalent to the conclusion of Theorem \ref{sumthm}. More
precisely we have

\begin{Thm}[{\cite[Twierdzenie 6.28]{phd}}]
Let $H$ be a Hilbert space and let $(A,B)$ be an irreducible pair of normal 
operators acting on $H$ satisfying relations \refeq{finalrel}. Assume that
$\spec{A}=\Ybar$. Then the following conditions are equivalent:
\begin{itemize}
\item[{\rm (1)}] the operator $A+B$ has a normal extension,
\item[{\rm (2)}] the operator $A\dplus B$ is normal,
\item[{\rm (3)}] $(A,B)\in\Dd_H$.
\end{itemize}
\end{Thm} 

With some more effort one can get rid of the assumption of irreducibility.
This theorem says that the spectral conditions imposed on pairs
$(A,B)\in\Dd_H$ are necessary for the sum $A\dplus B$ to have the same 
analytic properties as $A$ and $B$.

\Section{Affiliation relation}\label{affrel}

In this section we shall deal with the affiliation relation for
$\cst$-algebras investigated in \cite{unbo} and its relationship with the
special functions investigated in Section \ref{SpecFun}. We shall use the
notion of a $\cst$-algebra generated by a finite family of affiliated elements
as well as a $\cst$-algebra generated by a quantum family of affiliated 
elements. We refer to \cite{unbo} and \cite{gen} for a detailed exposition of
these topics.

A very simple but useful lemma presented below uses the interplay between
these concepts.

\begin{Lem}\label{afflem}
Let $H$ be a Hilbert space and let $\UA$ be a non degenerate $\cst$-subalgebra
of $\B{H}$. Let $\UC$ and $\UB$ be $\cst$-algebras and let 
$F\in\M{\UC\tens\UB}$ be a
quantum family of affiliated elements generating $\UB$. Let 
$\pi\in\Rep{\UB}{H}$
and let $R_1,\ldots,R_N$ be elements affiliated with $\UB$ such that
$R_1,\ldots,R_N$ generate $\UB$. Define $T_k=\pi(R_k)$ for $k=1,\ldots,N$. Then
\[
\Bigl(\:(\id_{\UC}\tens\pi)F\in\M{\UC\tens\UA}\:\Bigr)\Longleftrightarrow
\Bigl(\:T_k\aff\UA\text{ for }k=1,\ldots,N\:\Bigr)
\]
\end{Lem}

\proof
``$\Rightarrow$'': Since $F$ generates $\UB$, the condition that 
$(\id_{\UC}\tens\pi)F\in\M{\UC\tens\UA}$ implies that $\pi\in\Mor{\UB}{\UA}$.
Therefore for $k\in\{1,\ldots,N\}$ we have
\[
T_k=\pi(R_k)\aff\UA.
\]

``$\Leftarrow$'': As the operators $T_k=\pi(R_k)$  
are affiliated with $\UA$ and
$\{R_1,\ldots,R_N\}$ generate $\UB$, the representation $\pi$ is a
morphism from $\UB$ to $\UA$. Consequently
$(\id_{\UC}\tens\pi)\in\Mor{\UC\tens\UB}{\UC\tens\UA}$ and
\[
(\id_{\UC}\tens\pi)F\in\M{\UC\tens\UA}.
\]

\subsection{Generators of some $\cst$-algebras}

\begin{Prop}\label{Fgen}
Let $F$ be the element of
$\M{\Ci{\Ybar}\tens\Ci{\Ybar}}=\Cb{\Ybar\times\Ybar}$ given by
\[
F(\y,\y')=\Fq(\y\y').
\]
Then $F$ is a quantum family of elements generating $\Ci{\Ybar}$.
\end{Prop}

\proof
Since $F\in\M{\Ci{\Ybar}\tens\Ci{\Ybar}}$, the map
\[
\Ybar\ni\y\longmapsto F(\y,\cdot)\in\M{\Ci{\Ybar}}
\]
is strictly continuous (cf.~\cite[Section 2]{gen}).
In particular for any $\ph\in L^1(\Yq)$ we can consider the integral
\[
\Int_\Yq F(\y,\cdot)\ph(\y)\,dh(\y)\in\M{\Ci{\Ybar}}.
\]
Using the asymptotic behaviour
$\Fq(\y)\approx\alpha(q^{-1}\y)$ it is easy to show that the function
\begin{equation}\label{fph}
F_\ph\colon\Ybar\ni\y'\longmapsto\Int_\Yq F(\y,\y')\ph(\y)\,dh(\y)
\end{equation}
belongs to $\Ci{\Ybar}$. We shall show that the family of functions
$\{F(\y,\cdot)\}_{\y\in\Ybar}$ separates points of $\Ybar$.
Indeed: suppose that for some $\y_1,\y_2\in\Ybar$ we have
$F(\y,\y_1)=F(\y,\y_2)$
for all $\y\in\Ybar$. This means that 
\begin{equation}\label{rO}
\Fq(\y\y_1)=\Fq(\y\y_2)
\end{equation}
for all $\y\in\Ybar$. In particular for $\y=qq^{-it}$ ($t\in\RR$) we obtain
\[
\Fq^{q\y_1}(t)=\Fq(qq^{it}\y_1)=\Fq(qq^{it}\y_2)=\Fq^{q\y_2}(t)
\]
for all $t\in\RR$. Performing holomorphic continuation to $t=-i$ and using
\refeq{rownFq} we get
\[
(1+\y_1)\Fq(\y_1)=(1+\y_2)\Fq(\y_2),
\]
which by \refeq{rO} means that $\y_1=\y_2$.

It follows that the family of functions $\{F_\ph:\:\ph\in L^1(\Yq)\}$
also separates points of $\Ybar$ (e.g.~by considering integrable functions
approximating measures concentrated on single points of $\Ybar$). 
By Stone-Weierstrass theorem applied to one point compactification of
$\Ybar$ the $*$-algebra generated by functions of the form \refeq{fph} is
dense in $\Ci{\Ybar}$.

Now let $H$ be a Hilbert space and let 
$\pi\in\Rep{\Ci{\Ybar}}{H}$. Assume that for some non
degenerate $\cst$-subalgebra $\UB\subset\B{H}$ we have
\[
\bigl((\id\tens\pi)F\bigr)\in\M{\Ci{\Ybar}\tens\UB},
\]
i.e.~$\bigl((\id\tens\pi)F\bigr)$ is a strictly continuous function on $\Ybar$
with values in $\M{\UB}$. This function acts in the following way:
\[
\Ybar\ni\y\longmapsto\pi\bigl(F(\y,\cdot)\bigr)\in\M{\UB}.
\]
For $\ph\in L^1(\Yq)$ let us denote the functional
\[
\Ci{\Ybar}\ni f\longmapsto\Int_\Yq f(\y)\ph(\y)\,dh(\y)
\]
by $\omega_\ph$. Then
\[
\pi(F_\ph)=\pi\bigl((\omega_\ph\tens\id)F\bigr)
=(\omega_\ph\tens\id)\bigl((\id\tens\pi)F\bigr)\in\M{\UB}.
\]
Since the $*$-algebra generated by elements of $\{F_\ph:\:\ph\in L^1(\Yq)\}$ is
dense in $\Ci{\Ybar}$, we see that
\[
\pi\bigl(\Ci{\Ybar}\bigr)\subset\M{\UB}.
\]
It remains to show that the set $\pi\bigl(\Ci{\Ybar}\bigr)\UB$ is linearly dense
in $\UB$. We shall use the fact that $F$ is a unitary element of
$\M{\Ci{\Ybar}\tens\Ci{\Ybar}}$. It implies that $(\id\tens\pi)F$ is a unitary
element of $\M{\Ci{\Ybar}\tens\UB}$. In particular the set
\begin{equation}\label{zbg}
\left\{\bigl((\id\tens\pi)F\bigr)(f\tens x):\:f\in\Ci{\Ybar},\:
x\in\UB\right\}
\end{equation}
is linearly dense in $\Ci{\Ybar}\tens\UB$. Notice further that for any
functional $\omega\in\Ci{\Ybar}^*$ we have
\[
(\omega\tens\id)\bigl[\bigl((\id\tens\pi)F\bigr)(f\tens x)\bigr]=
\pi\bigl[(\omega\tens\id)\bigl(F(f\tens I)\bigr)\bigr]x.
\]
In particular for $\ph\in L^1(\Yq)$
\[
(\omega_{\ph}\tens\id)\bigl[\bigl((\id\tens\pi)F\bigr)(f\tens x)\bigr]=
\pi\bigl[(\omega_{f\ph}\tens\id)F\bigr]x\in\pi\bigl(\Ci{\Ybar}\bigr)\UB.
\]
For a given $y\in\UB$ consider $g\in\Ci{\Ybar}$ and $\ph\in L^1(\Yq)$ such that
$\omega_\ph(g)=1$. Choose a sequence $(q_n)_{n\in\NN}$ of linear
combinations of elements of \refeq{zbg} converging to $g\tens y$. Then the 
sequence of elements
\[
(\omega_\ph\tens\id)(q_n)\in\pi\bigl(\Ci{\Ybar}\bigr)\UB
\]
converges to $y$. This shows that $\pi\in\Mor{\Ci{\Ybar}}{\UB}$.
\qed

Using similar methods or appealing to the theory of multiplicative unitary 
operators (e.g.~\cite[Theorem 1.6(6)]{mu}) one can also prove:

\begin{Prop}
Let $F$ be the element of
$\M{\Ci{\Yq}\tens\Ci{\Yq}}=\Cb{\Yq\times\Yq}$ given by
\[
F(\y,\y')=\chi(\y,\y').
\]
Then $F$ is a quantum family of elements generating $\Ci{\Yq}$.
\end{Prop}

We can apply Lemma \ref{afflem} in the two situations described below.
\begin{itemize}
\item[({\em a}\/)] Let $H$ be a Hilbert space and let $T$ be a normal operator
acting on $H$ such that $\spec{T}\subset\Ybar$. Set $\UB=\UC=\Ci{\Ybar}$ and let
$F\in\M{\UC\tens\UB}=\Cb{\Ybar\times\Ybar}$ be given by 
\[
F(\y,\y')=\Fq(\y\y').
\] 
Then $F$ generates $\UB$.
The algebra $\UB$ is also generated by a single affiliated element $R_1$ given 
by $R_1(\y)=\y$. Further define the representation $\pi$ by $\pi(f)=f(T)$ for 
$f\in\UB$.
\item[({\em b}\/)] Let $H$ be a Hilbert space and let $T$ be a normal operator
acting on $H$ such that $\spec{T}\subset\Ybar$ and $\ker\{T\}=\{0\}$. 
Set $\UB=\UC=\Ci{\Yq}$ and let $F\in\M{\UC\tens\UB}=\Cb{\Yq\times\Yq}$ be given 
by
\[
F(\y,\y')=\chi(\y,\y').
\]
Then $F$ generates $\UB$. 
The algebra $\UB$ is also generated by two affiliated elements $R_1$
and $R_2$ given by $R_1(\y)=\y$, $R_2(\y)=\y^{-1}$. Further define a 
representation $\pi$ of $\UB$ by $\pi(f)=f(T)$ for $f\in\UB$.
\end{itemize}
Recall that for $\UC=\Ci{\Lambda}$ and any $\cst$-algebra $\UB$ the multiplier
algebra $\M{\UC\tens\UB}$ is canonically isomorphic to the algebra of strictly
continuous $\M{\UB}$-valued functions on $\Lambda$. 
Thus the two cases ({\em a}\/) and ({\em b}\/) yield the following theorems:

\begin{Thm}\label{af1}
Let $H$ be a Hilbert space, $T$ a normal operator acting on $H$ such that
$\spec{T}\subset\Ybar$ and let $\UA\subset B(H)$ be a non degenerate
$\cst$-subalgebra. Then the following conditions are equivalent:
\begin{itemize}
\item[{\rm (1)}] for any $\y\in\Ybar$
the unitary operator $\Fq(\y T)$ belongs to $\M{\UA}$ and the map
\[
\Ybar\ni\y\longmapsto\Fq(\y T)\in\M{\UA}
\]
is strictly continuous;
\item[{\rm (2)}] the operator $T$ is affiliated with $\UA$.
\end{itemize}
\end{Thm}

\begin{Thm}\label{af2}
Let $H$ be a Hilbert space, $T$ a normal operator acting on $H$ such that
$\spec{T}\subset\Ybar$ and\/ $\ker{T}=\{0\}$. Let $\UA\subset B(H)$ be a non 
degenerate $\cst$-subalgebra. Then the following conditions are equivalent:
\begin{itemize}
\item[{\rm (1)}] for any $\y\in\Yq$
the unitary operator $\chi(\y,T)$ belongs to $\M{\UA}$ and the map
\[
\Yq\ni\y\longmapsto\chi(\y,T)\in\M{\UA}
\]
is strictly continuous;
\item[{\rm (2)}] operators $T$ and $T^{-1}$ are affiliated with $\UA$.
\end{itemize}
\end{Thm}

\subsection{Algebraic consequences}

\begin{Thm}
Let $H$ be a Hilbert space, $\UA\subset\B{H}$ a non degenerate $\cst$-subalgebra
and let $(A,B)\in\Dd_H$ be such that $A,B\aff\UA$. Then
\begin{itemize}
\item[{\rm (1)}] the operator $A\dplus B$ is affiliated with $\UA$,
\item[{\rm (2)}] the operator $BA$ is affiliated with $\UA$.
\end{itemize} 
\end{Thm}

\proof
{\sc Ad (1).} By Statement (1) of Corollary \ref{corAB} for any 
$\y\in\Yq$ we have $(\y A,\y B)\in\Dd_H$ and thus by \refeq{expon}
\[
\Fq\bigl(\y(A\dplus B)\bigr)=\Fq(\y B)\Fq(\y A).
\]
Now the result follows by Theorem \ref{af1}.

{\sc Ad (2).} We apply the same method as in the proof of (1) and use
\refeq{F5} and Theorem \ref{af2}.
\qed

We end this section with the following useful proposition.  

\begin{Prop}\label{X}
Let $X$ be a closed subset of\/ $\CC\bez\{0\}$. Let $H$ be a
Hilbert space and let $T_1,T_2$ be strongly commuting normal operators acting 
on $H$ such that\/ $\spec{T_1},\:\spec{T_2}\subset\overline{X}$ and\/
$\ker{T_1}=\ker{T_2}=\{0\}$. Let $\UA$ be a non degenerate $\cst$-subalgebra of
$\B{H}$ and assume that $T_1,T_2,T_1^{-1}$ and $T_2^{-1}$ are affiliated with
$\UA$. Then for any $f\in\Cb{X\times X}$ we have $f(T_1,T_2)\in\M{\UA}$.
\end{Prop}

The proof of this proposition is identical to that of
\cite[Proposition 5.3]{azb} (cf.~also \cite{phd}). We have formulated it in a
way which gets rid of unnecessary assumptions about the shape of $X$. 

We shall use this proposition with $X=\Yq$.

\Section{Multiplicative unitary and its properties}\label{MU}

Most results of this section are direct analogues of those presented in
\cite[Sections 3--7]{azb}. Also most proofs are identical. Therefore we
shall omit some of the proofs. It needs to be stressed that this analogy stems
from the fact that the formulae arising in the study of commutation
relations described in Subsection \ref{komrel} are in most cases identical to
those found in \cite{azb}. However the meaning of these formulae is different,
as the commutation relations discussed in both papers are different. We shall
include the proofs of the results relying on the theorems proved in this
paper.

\subsection{The quantum group space}

As in \cite{azb} we shall begin the construction of our quantum ``$az+b$''
group with the definition of the operator domain $\Gg$ playing the role of the
quantum space of our quantum group. 

Let $H$ be a Hilbert space. By $\Gg_H$ we shall denote the set of pairs $(a,b)$
of closed operators on $H$ satisfying
\begin{enumerate}
\item $a$ and $b$ are normal,
\item $\spec{a},:\spec{b}\subset\Ybar$,
\item $\ker{a}=\{0\}$,
\item $\chi(\y,a)b\chi(\y,a)^*=\y b$ for all $\y\in\Yq$.
\end{enumerate}
It follows that $a$ preserves the decomposition
$H=\ker{b}\oplus(\ker{b})^{\ort}$
and $(a,b)\in\Gg_H$ if and only if $a_0=\bigl.a\bigr|_{\ker{b}}$ satisfies
$\ker{a_0}=\{0\}$, $\spec{a_0}\subset\Ybar$ and the pair
$(A,B)=(\bigl.a\bigr|_{(\ker{b})^\ort},\bigl.b\bigr|_{(\ker{b})^\ort})$ belongs 
to $\Dd_{(\ker{b})^\ort}$.

The operator domain $\Gg$ will serve as the quantum space of our quantum group
in the sense that to each representation of the algebra of continuous functions 
vanishing at infinity on the quantum group on a Hilbert space $H$ there will 
correspond a unique pair $(a,b)\in\Gg_H$ (cf.~Subsection \ref{alg}). 

It should be pointed out that the definition of $\Gg$ is not an essential
ingredient of the construction of the new quantum ``$az+b$'' groups. In
practice we can always choose a faithful representation with $(a,b)\in\Dd_H$
(cf.~remarks after Proposition \ref{univ}).

\begin{Lem}\label{ABKOM}
Let $H$ and $K$ be Hilbert spaces and let $(a,b)\in\Gg_H$ and
$(\ahat,\bhat)\in\Gg_K$. Then there is the following relation between
operators on $K\tens H$:
\[
\chi(\ahat\tens I,I\tens a)(\bhat\tens I)=(\bhat\tens a)
\chi(\ahat\tens I,I\tens a).
\]
\end{Lem}

The assertion of Lemma \ref{ABKOM} follows for example from 
\cite[Formula 2.6]{gen} (cf.~also \cite{azb,axb}).

\subsection{Multiplicativity}

\begin{Prop}
Let $H$ be a Hilbert space and let $(a,b)\in\Gg_H$ be such that
$\ker{b}=\{0\}$. Then the unitary operator
\begin{equation}\label{W}
W=\Fq(b^{-1}a\tens b)\chi(b^{-1}\tens I,I\tens a)
\end{equation}
on $H\tens H$ satisfies
\begin{equation}\label{Wkomult}
\begin{array}{r@{\;=\;}l@{\smallskip}}
W(a\tens I)W^*&a\tens a,\\
W(b\tens I)W^*&a\tens b\dplus b\tens a. 
\end{array}
\end{equation}
\end{Prop}

\proof
It is easy to see that $(b^{-1},a)\in\Gg_H$. By Lemma \ref{ABKOM}
\[
\chi(b^{-1}\tens I,I\tens a)(a\tens I)\chi(b^{-1}\tens I,I\tens a)^*
=a\tens a.
\]
Also since $a\tens a$ strongly commutes with $b^{-1}a\tens b$, we have
\[
\begin{array}{r@{\;=\;}l@{\smallskip}}
W(a\tens I)W^*&\Fq(b^{-1}a\tens b)\chi(b^{-1}a\tens I,I\tens a)(a\tens I)
\chi(b^{-1}\tens I,I\tens a)^*\Fq(b^{-1}a\tens b)^*\\
&\Fq(b^{-1}\tens b)(a\tens a)\Fq(b^{-1}a\tens b)^*\\
&a\tens a.
\end{array}
\]

For the second formula of \refeq{Wkomult} notice that the pair
$(A,B)=(a\tens b,b\tens I)\in\Dd_{H\tens H}$. Also $b\tens I$ commutes with 
$\chi(b^{-1}\tens I,I\tens a)$. Therefore by \refeq{suma2}
\[
\begin{array}{r@{\;=\;}l@{\smallskip}}
W(b\tens I)W^*&\Fq(b^{-1}a\tens b)\chi(b^{-1}\tens I,I\tens a)(b\tens I)
\chi(b^{-1}\tens I,I\tens a)^*\Fq(b^{-1}a\tens b)^*\\
&\Fq(b^{-1}a\tens b)(b\tens I)\Fq(b^{-1}a\tens b)^*\\
&\Fq(B^{-1}A)B\Fq(B^{-1}A)^*=A\dplus B\\
&a\tens b\dplus b\tens I.
\end{array}
\]
\qed

\begin{Prop}\label{adapProp}
Let $H$ and $K$ be Hilbert spaces and let $(a,b)\in\Gg_H$ and
$(\ahat,\bhat)\in\Gg_K$. Assume that $\ker{b}=\{0\}$. Then the operator $W$
defined by \refeq{W} and 
\begin{equation}\label{V}
V=\Fq(\bhat\tens b)\chi(\ahat\tens I,I\tens a)
\end{equation} 
satisfy
\begin{equation}\label{adap}
W_{23}V_{12}=V_{12}V_{13}W_{23}
\end{equation}
on $K\tens H\tens H$.
\end{Prop}

\proof
Using \refeq{Wkomult} we obtain
\[
\begin{array}{r@{\;}c@{\;}l@{\smallskip}}
W_{23}V_{12}W_{23}^*&=&
(I\tens W)\bigl(\Fq(\bhat\tens I)\tens I\bigr)(I\tens W)^*\\
&&\times(I\tens W)\bigl(\chi(\ahat\tens I,I\tens a)\tens I\bigr)(I\tens W)^*\\
&=&\Fq\bigl(\bhat\tens(a\tens b\dplus b\tens I)\bigr)
\chi(\ahat\tens I\tens I,I\tens a\tens a).
\end{array}
\]
Also
\[
\chi(\ahat\tens I\tens I,I\tens a\tens a)=
\chi(\ahat\tens I\tens I,I\tens a\tens I)
\chi(\ahat\tens I\tens I,I\tens I\tens a),
\]
since $\chi$ is a bicharacter on $\Yq$.

We shall consider two cases: $\bhat=0$ and $\ker{\bhat}=\{0\}$. Assume first
that $\bhat=0$. Then $V=\chi(\ahat\tens I,I\tens a)$ and
\[
\begin{array}{r@{\;=\;}l@{\smallskip}}
W_{23}V_{12}W_{23}^*&\chi(\ahat\tens I\tens I,I\tens a\tens a)\\
&\chi(\ahat\tens I\tens I,I\tens a\tens I)
\chi(\ahat\tens I\tens I,I\tens I\tens a)\\
&V_{12}V_{13}.
\end{array}
\]
In the other case $(\ahat,\bhat)\in\Dd_K$ and therefore
\begin{equation}\label{noweRS}
(A,B)=(\bhat\tens a\tens b,\bhat\tens b\tens I)\in\Dd_{K\tens H\tens H}.
\end{equation}
Lemma \ref{ABKOM} says that
$(\bhat\tens a)\chi(\ahat\tens I,I\tens a)
=\chi(\ahat\tens I,I\tens a)(\bhat\tens I)$,
and consequently
\begin{equation}\label{potrzeb}
(\bhat\tens a\tens b)\chi(\ahat\tens I\tens I,I\tens a\tens I)
=\chi(\ahat\tens I\tens I,I\tens a\tens I)(\bhat\tens I\tens b)
\end{equation}
Now using \refeq{expon} with the pair \refeq{noweRS} we get:
\[
\Fq\bigl(\bhat\tens(a\tens b\dplus b\tens I)\bigr)
=\Fq(\bhat\tens a\tens b\dplus\bhat\tens b\tens I)
=\Fq(\bhat\tens b\tens I)\Fq(\bhat\tens a\tens b).
\]
Therefore
\[
\begin{array}{r@{\;=\;}l@{\smallskip}}
W_{23}V_{12}W_{23}^*&\Fq\bigl(\bhat\tens(a\tens b\dplus b\tens I)\bigr)
\chi(\ahat\tens I\tens I,I\tens a\tens a)\\
&\Fq\bigl(\bhat\tens(a\tens b\dplus b\tens I)\bigr)
\chi(\ahat\tens I\tens I,I\tens a\tens I)
\chi(\ahat\tens I\tens I,I\tens I\tens a)\\
&\Fq(\bhat\tens b\tens I)\Fq(\bhat\tens a\tens b)
\chi(\ahat\tens I\tens I,I\tens a\tens I)
\chi(\ahat\tens I\tens I,I\tens I\tens a)\\
&\Fq(\bhat\tens b\tens I)
\chi(\ahat\tens I\tens I,I\tens a\tens I)
\Fq(\bhat\tens I\tens b)
\chi(\ahat\tens I\tens I,I\tens I\tens a)\\
&V_{12}V_{13}.
\end{array}
\]
where in the second last equality we used \refeq{potrzeb}.

Now for general $\bhat$ we split the Hilbert space $K\tens H\tens H$ into a 
direct sum of $(\ker{\bhat})\tens H\tens H$ and 
$(\ker{\bhat})^\ort\tens H\tens H$ and derive \refeq{adap} separately for each
summand.
\qed

Setting $K=H$ and $(\ahat,\bhat)=(b^{-1},b^{-1}a)$ in Proposition 
\ref{adapProp} we obtain $V=W$ and thus we prove Statement (1) of the
following corollary:

\begin{Cor}
Let $H$ be a Hilbert space and let $(a,b)\in\Gg_H$ be such that
$\ker{b}=\{0\}$. Then
\begin{itemize}
\item[{\rm (1)}] The operator $W$ defined by \refeq{W} is a multiplicative
unitary,
\item[{\rm (2)}] for any Hilbert space $K$ and any $(\ahat,\bhat)\in\Gg_H$ the
operator $V$ defined by \refeq{V} is a unitary adapted to $W$.
\end{itemize}
\end{Cor}

\subsection{Modularity}

In this subsection we shall show that the multiplicative unitary operator
given by \refeq{W} is {\em modular}\/ (cf.~\cite[Definition 2.1]{modmu}).
In what follows we shall need the {\em partial transposition formula}\/
(cf.~\cite[Formula 3.15]{azb} or \cite[Lemat 7.7]{phd}). More precisely we 
shall use the
following fact: let $H$ and $K$ be Hilbert spaces and let $a$ and $\ahat$ be
normal operators acing on $H$ and $K$ respectively. Then for any bounded Borel
function $f$ on $\spec{\ahat}\times\spec{a}$ and all $z,x\in K$, $u,y\in H$ we
have
\begin{equation}\label{partrwz}
\its{\zbar\tens u}{f(\ahat^{\top}\tens I,I\tens a)}{\xbar\tens y}=
\its{x\tens u}{f(\ahat\tens I,I\tens a)}{z\tens y}.
\end{equation}

\begin{Prop}\label{zwiazek}
Let $H$ and $K$ be Hilbert spaces and let $(a,b)\in\Gg_H$,
$(\ahat,\bhat)\in\Gg_K$. Assume also that $\ker{b}=\{0\}$ and
$\ker{\bhat}=\{0\}$. Further denote $Q=|a|$ and let $x,z\in K$, $u\in D(Q)$
and $y\in D(Q^{-1})$. Define functions $\ph,\psi\colon\Yq\to\CC$ as
\begin{equation}\label{phipsi}
\begin{array}{r@{\;=\;}l@{\medskip}}
\ph(\y)&\its{x\tens u}{\chi(\bhat\tens b,\y)
\chi(\ahat\tens I,I\tens a)}{z\tens y},\\
\psi(\y)&\its{\zbar\tens Qu}{\chi(-\bhat^{\top}\tens qa^{-1}b,\y)
\chi(\ahat^{\top}\tens I,I\tens a)}{\xbar\tens Q^{-1}y}.
\end{array}
\end{equation}
Then
\begin{itemize}
\item[{\rm (1)}] $\psi(\y)=\overline{\alpha(\y)}|\y|\chi(-1,\y)\ph(\y)$,
\item[{\rm (2)}] if $x\in D(\bhat^{\pm1})$, $u\in D(b^{\pm1}\comp Q^{\pm2})$
and $y\in D(Q^{\pm2})$ for all possible combinations of signs then $\ph$ and
$\psi$ belong to the space $\Sch(\Yq)$.
\end{itemize}
\end{Prop}

\proof
{\sc Ad (1).} Since $(a,b)\in\Dd_H$ and $Q=|a|$, we have
\[
(I\tens Q^{-it})(-\bhat^{\top}\tens qa^{-1}b)(I\tens Q^{it})
=q^{-it}(-\bhat^{\top}\tens qa^{-1}b)
\]
(cf.~\refeq{finalrel}). Applying to both sides of this equation the function
$\y'\mapsto\chi(\y',\y)$ we obtain
\begin{equation}\label{komQ}
\begin{array}{r@{}l@{\smallskip}}
(I\tens Q^{-it})&\chi(-\bhat^{\top}\tens qa^{-1}b,\y)(I\tens Q^{it})\\
&=\chi(q^{-it}(-\bhat^{\top}\tens qa^{-1}b),\y)\\
&=\chi(q^{-it},\y)\chi(-\bhat^{\top}\tens qa^{-1}b,\y)\\
&=|\y|^{-it}\chi(-\bhat^{\top}\tens qa^{-1}b,\y).
\end{array}
\end{equation}
Now $Q$ commutes with $a$ and consequently $(I\tens Q^{it})$
commutes with $\chi(\ahat^{\top}\tens I,I\tens a)$. Therefore
\[
\begin{array}{l@{\medskip}}
\its{\zbar\tens Q^{it}u}{\chi(-\bhat^{\top}\tens qa^{-1}b,\y)
\chi(\ahat^{\top}\tens I,I\tens a)}{\xbar\tens Q^{it}y}\\
\qquad=\its{\zbar\tens u}{(I\tens Q^{-it})\chi(-\bhat^{\top}\tens qa^{-1}b,\y)
\chi(\ahat^{\top}\tens I,I\tens a)(I\tens Q^{it})}{\xbar\tens y}\\
\qquad=\its{\zbar\tens u}{(I\tens Q^{-it})\chi(-\bhat^{\top}\tens qa^{-1}b,\y)
(I\tens Q^{it})\chi(\ahat^{\top}\tens I,I\tens a)}{\xbar\tens y}\\
\qquad=|\y|^{-it}\its{\zbar\tens u}{\chi(-\bhat^{\top}\tens qa^{-1}b,\y)
\chi(\ahat^{\top}\tens I,I\tens a)}{\xbar\tens y}.
\end{array}
\]
Performing holomorphic continuation to $t=i$ we get
\[
\begin{array}{l@{\medskip}}
\its{\zbar\tens Qu}{\chi(-\bhat^{\top}\tens qa^{-1}b,\y)
\chi(\ahat^{\top}\tens I,I\tens a)}{\xbar\tens Q^{-1}y}\\
\qquad=|\y|\its{\zbar\tens u}{\chi(-\bhat^{\top}\tens qa^{-1}b,\y)
\chi(\ahat^{\top}\tens I,I\tens a)}{\xbar\tens y}.
\end{array}
\]
Notice that
\[
\begin{array}{r@{\;=\;}l@{\smallskip}}
\chi(-\bhat^{\top}\tens qa^{-1}b,\y)
&\chi(-1,\y)\chi(\bhat^{\top}\tens qa^{-1}b,\y)\\
&\chi(-1,\y)\chi(\bhat^{\top}\tens I,\y)\chi(I\tens qa^{-1}b,\y)\\
&\chi(-1,\y)\chi(\bhat,\y)^{\top}\tens\chi(qa^{-1}b,\y).
\end{array}
\]
Therefore we can continue our computation in the following way:
\[
\begin{array}{l@{\medskip}}
\its{\zbar\tens Qu}{\chi(-\bhat^{\top}\tens qa^{-1}b,\y)
\chi(\ahat^{\top}\tens I,I\tens a)}{\xbar\tens Q^{-1}y}\\
\qquad=|\y|\chi(-1,\y)\its{\zbar\tens u}{
\left(\chi(\bhat,\y)^{\top}\tens\chi(qa^{-1}b,\y)\right)
\chi(\ahat^{\top}\tens I,I\tens a)}{\xbar\tens y}\\
\qquad=|\y|\chi(-1,\y)\its{\overline{\chi(\bhat,\y)z}\tens \chi(qa^{-1}b,\y)^*u}
{\chi(\ahat^{\top}\tens I,I\tens a)}{\xbar\tens y}\\
\qquad=|\y|\chi(-1,\y)\its{x\tens \chi(qa^{-1}b,\y)^*u}
{\chi(\ahat\tens I,I\tens a)}{\chi(\bhat,\y)z\tens y},
\end{array}
\]
where in the last equality we used \refeq{partrwz}.
Then we have
\[
\begin{array}{l@{\medskip}}
\its{\zbar\tens Qu}{\chi(-\bhat^{\top}\tens qa^{-1}b,\y)
\chi(\ahat^{\top}\tens I,I\tens a)}{\xbar\tens Q^{-1}y}\\
\qquad=|\y|\chi(-1,\y)\its{x\tens u}{\bigl(I\tens\chi(qa^{-1}b,\y)\bigr)
\chi(\ahat\tens I,I\tens a)\bigl(\chi(\bhat,\y)\tens I\bigr)}{z\tens y}\\
\qquad=|\y|\chi(-1,\y)\its{x\tens u}{\bigl(I\tens\chi(qa^{-1}b,\y)\bigr)
\chi(\ahat\tens I,I\tens a)\chi(\bhat\tens I,\y)}{z\tens y}.
\end{array}
\]
By Lemma \ref{ABKOM} this last expression equals
\[
|\y|\chi(-1,\y)\its{x\tens u}{\bigl(I\tens\chi(qa^{-1}b,\y)\bigr)
\chi(\bhat\tens a,\y)\chi(\ahat\tens I,I\tens a)}{z\tens y}.
\]
Let us put $A=\bhat\tens a$, $B=I\tens qa^{-1}b$. It is easy to see that
$(A,B)\in\Dd_{K\tens H}$. By Statement (4) of Corollary \ref{corAB}
\[
\begin{array}{r@{\;=\;}l@{\smallskip}}
\chi(I\tens qa^{-1}b,\y)\chi(\bhat\tens a,\y)
&\chi(B,\y)\chi(\y,A)\\
&\overline{\alpha(\y)}\chi(qBA,\y)\\
&\overline{\alpha(\y)}\chi(\bhat\tens b,\y).
\end{array}
\]
With this information we obtain
\[
\begin{array}{l@{\medskip}}
\its{\zbar\tens Qu}{\chi(-\bhat^{\top}\tens qa^{-1}b,\y)
\chi(\ahat^{\top}\tens I,I\tens a)}{\xbar\tens Q^{-1}y}\\
\qquad=\overline{\alpha(\y)}|\y|\chi(-1,\y)
\its{x\tens u}{\chi(\bhat\tens b,\y)\chi(\ahat\tens I,I\tens a)}{z\tens y},
\end{array}
\]
which proves Statement (1).

{\sc Ad (2).}
In a similar way to the derivation of \refeq{komQ} we find that
\[
(I\tens Q^{-it})\chi(\bhat\tens b,\y)(I\tens Q^{it})
=|\y|^{-it}\chi(\bhat\tens b,\y)
\]
for all $\y\in\Yq$ and all $t\in\RR$. Therefore
\[
\begin{array}{l@{\smallskip}}
\its{x\tens Q^{it}u}{\chi(\bhat\tens b,\y)\chi(\ahat\tens I,I\tens a)}
{z\tens Q^{it}y}\\
\qquad=|\y|^{-it}\its{x\tens u}{\chi(\bhat\tens b,\y)
\chi(\ahat\tens I,I\tens a)}
{z\tens y}=|\y|^{-it}\ph(\y).
\end{array}
\]
Performing holomorphic continuation to the point $t=\pm2i$ we obtain
\[
\its{x\tens Q^{\pm2}u}{\chi(\bhat\tens b,\y)\chi(\ahat\tens I,I\tens a)}
{z\tens Q^{\mp2}y}=|\y|^{\pm2}\ph(\y).
\]
With $\y=q^kq^{it}$ this means that
\[
e^{\mp2t\iim{\rh^{-1}}}\ph\bigl(q^kq^{it}\bigr)
=\its{x\tens Q^{\pm2}u}{\PHase{\bhat\tens b}|\bhat\tens b|^{it}
\chi(\ahat\tens I,I\tens a)}{z\tens Q^{\mp2}y}.
\]
It follows from the assumptions about $x,y$ and $u$ that $x\tens
Q^{\pm2}u\in\D{\bhat\tens b}$. Consequently the functions
\[
t\longmapsto e^{\mp2t\iim{\rh^{-1}}}\ph\bigl(q^kq^{it}\bigr)
\]
have a holomorphic continuation to $\{z\in\CC:\:-1<\im{z}<1\}$ and this
continuation is bounded in the strip. It is an easy exercise (cf.~\cite[Lemat
7.9]{phd}) that if a function $u$ on $\RR$ has the property that the functions
$t\mapsto e^{\pm t}u(t)$ extend to bounded holomorphic functions on the strip
$\{z\in\CC:\:-1<\im{z}<1\}$ then $u\in\Sch(\RR)$. In particular the functions 
$t\mapsto\ph(q^kq^{it})$ belong to $\Sch(\RR)$, i.e.~$\ph\in\Sch(\Yq)$.

By Statement (1) we know that the functions $t\longmapsto\psi(q^kq^{it})$ are
smooth and that $t\mapsto e^{\pm t\iim{\rh^{-1}}}\psi(q^kq^{it})$ have
extensions to bounded holomorphic functions in the strip 
$\{z\in\CC:-1<\im{z}<1\}$. Therefore $\psi\in\Sch(\Yq)$.
\qed

In the same way as \cite[Proposition 3.5]{azb} we get the following:

\begin{Prop}\label{fig}
Let $H$, $K$, $(a,b)$, $(\ahat,\bhat)$, $Q$ and $x,z,u,y$ be as in 
Proposition \ref{zwiazek}. Let $f$ and $g$ be bounded Borel functions on $\Yq$.
Denote by $\widehat{f}$ and $\widehat{g}$ the inverse Fourier transforms of
$f$ and $g$ (we treat $f$ and $g$ as tempered distributions on $\Yq$):
\begin{equation}\label{Ftr}
\begin{array}{r@{\;=\;}l@{\smallskip}}
\displaystyle{f(\y)}&
\displaystyle{\Int_{\Yq}\widehat{f}(\y')\chi(\y,\y')\,d\Mu(\y'),}\\
\displaystyle{g(\y)}&
\displaystyle{\Int_{\Yq}\widehat{g}(\y')\chi(\y,\y')\,d\Mu(\y').}
\end{array}
\end{equation}
Suppose that for almost all $\y\in\Yq$ we have
\begin{equation}\label{zw2}
\widehat{f}(\y)=\overline{\alpha(\y)}|\y|\chi(-1,\y)\widehat{g}(\y).
\end{equation}
Then 
\begin{equation}\label{cel}
\begin{array}{l@{\medskip}}
\its{x\tens u}{f(\bhat\tens b,\y)\chi(\ahat\tens I,I\tens a)}{z\tens y}\\
\qquad=\its{\zbar\tens Qu}{g(-\bhat^{\top}\tens qa^{-1}b,\y)
\chi(\ahat^{\top}\tens I,I\tens a)}{\xbar\tens Q^{-1}y}.
\end{array}
\end{equation}
\end{Prop}

\begin{Cor}
Let $H$ and $K$ be Hilbert spaces and let $(a,b)\in\Gg_H$,
$(\ahat,\bhat)\in\Gg_K$. Assume that\/ $\ker{b}=\{0\}$ and\/ 
$\ker{\bhat}=\{0\}$.
Let $V$ be the operator introduced by \refeq{V} and
define
\[
\Vtil=\Fq(-\bhat^\top\tens qa^{-1}b)^*\chi(\ahat^\top\tens I,I\tens a).
\]
Then 
\begin{itemize}
\item[{\rm (1)}] for all $x,z\in K$, $y\in D(Q^{-1})$ and $u\in D(Q)$ we have
\begin{equation}\label{modu}
\its{x\tens u}{V}{z\tens y}=\its{\zbar\tens Qu}{\Vtil}{\xbar\tens Q^{-1}y};
\end{equation}
\item[{\rm (2)}] the operator $W$ introduced by \refeq{W} is a modular
multiplicative unitary.
\end{itemize}
\end{Cor}

\proof
{\sc Ad (1).} We shall use Proposition \ref{fig} with $f=\Fq$ and
$g=\overline{\Fq}$. We can use it because of Corollary \ref{wazne}. The result
is exactly \refeq{modu}.

{\sc Ad (2).} Put $K=H$ and $(\ahat,\bhat)=(b^{-1},b^{-1}a)$. Then by
Statement (1) we have
\[
\its{x\tens u}{W}{z\tens y}=\its{\zbar\tens Qu}{\Wtil}{\xbar\tens Q^{-1}y}
\]
for all $x,z\in K$, $y\in D(Q^{-1})$ and $u\in D(Q)$. Let $\Qhat=|b|$. Then it
is easy to verify that $\Qhat\tens Q$ strongly commutes with 
$b^{-1}a\tens b$, $b^{-1}\tens I$ and $I\tens a$. Therefore
\[
W(\Qhat\tens Q)W^*=\Qhat\tens Q.
\] 
This way we have verified that $W$ satisfies all conditions listed in
\cite[Definition 2.1]{modmu}.
\qed

\Section{The quantum ``$az+b$'' group for new values of $q$}\label{NewQG}

\subsection{The $\cst$-algebra}\label{alg}

Let us describe the $\cst$-algebra which will turn out to be the algebra of
continuous functions vanishing at infinity on the quantum ``$az+b$'' group.
For $\y\in\Yq$ and $f\in\Ci{\Ybar}$ let
\[
\bigl(\beta_\y f\bigr)(\y')=f(\y'\y)
\]
for all $\y'\in\Ybar$. Then $\Yq\ni\y\mapsto\beta_\y\in\Aut{\Ci{\Ybar}}$ is a
strongly continuous action and $(\Ci{\Ybar},\Yq,\beta)$ is a $\cst$-dynamical
system. Let $\UB$ be the corresponding $\cst$-crossed product. As the
canonical embedding $\Ci{\Ybar}\hookrightarrow\M{\UB}$ is a morphism from
$\Ci{\Ybar}$ to $\UB$, any element affiliated with $\Ci{\Ybar}$ can be treated
as an element affiliated with $\UB$. Let $b$ be the element affiliated with
$\UB$ arising from the continuous function $\Ybar\ni\y\mapsto\y\in\CC$. Let
$(U_\y)_{\y\in\Yq}$ be the strictly continuous family of unitary elements of 
$\M{\UB}$ implementing the action $\beta$: $U_\y fU_{\y}^*=\beta_\y f$ for
$f\in\Ci{\Ybar}\subset\M{\UB}$. Let us represent $\UB$ faithfully on a Hilbert
space $H$. Then $(U_\y)_{\y\in\Yq}$ is a group of unitary operators acting on
$H$. By SNAG theorem there is a normal operator $a$ acting on $H$ such that
$\spec{a}\subset\Ybar$, $\ker{a}=\{0\}$ and $U_\y=\chi(\y,a)$ for all
$\y\in\Yq$ (remark that this is where we are using the fact that
$\widehat{\Yq}=\Yq$). Now since all operators $U_\y$ are in $\M{\UB}$ and the
map $\Yq\ni\y\mapsto U_\y\in\M{\UB}$ is strictly continuous, by Theorem
\ref{af2} the operators $a$ and $a^{-1}$ are affiliated with $\UB$. It is easy
to check that for any $\pi\in\Rep{\UB}{H}$ we have $(\pi(a),\pi(b))\in\Gg_H$.
It is also known that 
\begin{equation}\label{dense}
\bigl\{g(a)f(b):\:g\in\Ci{\Yq},\:f\in\Ci{\Ybar}\bigr\}
\end{equation}
is a linearly dense subset of $\UB$.

Using the same technique as in the proof of 
\cite[Propositions 4.1 and 4.2]{azb} we get
the following:

\begin{Prop}\label{univ}
Let $\UB$, $a$ and $b$ be the $\cst$-algebra and two affiliated elements
described in this subsection. Then
\begin{itemize}
\item[{\rm (1)}] the $\cst$-algebra $\UB$ is generated by the three affiliated
elements $a$, $a^{-1}$ and $b$;
\item[{\rm (2)}] for any Hilbert space $H$ and any $(a_0,b_0)\in\Gg_H$ there
exists a unique $\pi\in\Rep{\UB}{H}$ such that $a_0=\pi(a)$ and $b_0=\pi(b)$.
If\/ $\UA$ is a non degenerate $\cst$-subalgebra of $\B{H}$ and
$a_0,a_0^{-1},b\aff\UA$ then $\pi\in\Mor{\UB}{\UA}$. 
\end{itemize}
\end{Prop}

Assume now that $\UB$ is faithfully represented in a Hilbert space $H$. From
the commutation relations between $a$ and $b$ and the fact that \refeq{dense}
is linearly dense in $\UB$ we see that $\ker{b}$ is an
invariant subspace of $H$ for the action of $\UB$. We can therefore restrict
our representation to $(\ker{b})^\ort$ or, equivalently, assume that
$\ker{b}=\{0\}$. Denote $\ahat=b^{-1}$ and $\bhat=b^{-1}a$ then 
\[
W=\Fq(\bhat\tens b)\chi(\ahat\tens I,I\tens a)
\]
coincides with \refeq{W}. Therefore it is a modular multiplicative unitary.
The operators $Q$, $\Qhat$ and $\Wtil$ related to $W$ via 
\cite[Definition 2.1]{azb} are given by
\[
\begin{array}{c@{\smallskip}}
\Qhat=|b|,\qquad Q=|a|,\\
\Wtil=\Fq(-\bhat^\top\tens qa^{-1}b)^*\chi(\ahat^\top\tens I,I\tens a).
\end{array}
\]
Let
\[
\UA=\bigl\{(\omega\tens\id)W:
\:\omega\in\B{H}_*\bigr\}^{\text{\rm\tiny norm closure}}
\!\!\!\!\!\!\!\!\!\!\!\!\!\!\!\!\!\!\!\!\!\!.
\]
By the theory developed in \cite{mu,modmu} $\UA$ is a nondegenerate
$\cst$-subalgebra of $\B{H}$.

The proof of the following proposition is virtually identical to that
presented in \cite[Section 6]{azb}. We give it here because it relies on the
results about special functions and commutation relations obtained in Sections 
\ref{SpecFun} and \ref{opeq}.

\begin{Prop}
The $\cst$-algebras $\UA$ and $\UB$ are equal as subsets of $\B{H}$.
\end{Prop}

\proof
The operators $\ahat\tens I$, $I\tens a$, $I\tens a^{-1}$ and $\bhat\tens b$
are affiliated with $\K{H}\tens\UB$. Therefore by Proposition \ref{X} and
Theorem \ref{af1} we have
\[
\chi(\ahat\tens I,I\tens a),\:\Fq(\bhat\tens b)\in\M{\K{H}\tens\UB}
\] 
and it follows that $W\in\M{\K{H}\tens\UB}$. By the definition of $\UA$ we
conclude that $\UA\subset\M{\UB}$. In particular $\UA\UB\subset\UB$. Using the
same technique as in the proof of Proposition \ref{Fgen} (cf.~also
\cite[Section 4]{mu} and \cite[Section 6]{azb}) one can show that $\UA\UB$ is
dense in $\UB$.

Now let us show that the elements $a$, $a^{-1}$ and $b$ are affiliated with
$\UA$. For any $\y\in\Ybar$ let
\[
V(\y)=\Fq(\y\bhat\tens b)\chi(\ahat\tens I,I\tens a).
\]
It is easy to verify with Theorem \ref{af1} that 
$\bigl(V(\y)\bigr)_{\y\in\Ybar}$ is a 
strictly
continuous family of unitary elements of $\M{\K{H}\tens\K{H}}$. Clearly if
$V(\y)_{12}=V(\y)\tens I$ then $\bigl(V(\y)_{12}\bigr)_{\y\in\Ybar}$ 
is a strictly continuous family of unitary elements
$\M{\K{H}\tens\K{H}\tens\UA}$. By Proposition \ref{adapProp}
\[
V(\y)_{13}=V(\y)_{12}W_{23}V(\y)_{12}W_{23}^*.
\]
Since $W\in\M{\K{H}\tens A}$
the family $\bigl(V(\y)_{13}\bigr)_{\y\in\Ybar}$ is a strictly continuous 
family of
elements of $\M{\K{H}\tens\K{H}\tens\UA}$. It follows that 
$\bigl(V(\y)\bigr)_{\y\in\Ybar}$ is a strictly continuous family of unitaries
in $\M{\UA}$. This implies that 
$\Fq(\y\bhat\tens b)=V(\y)V(0)^*\in\M{\K{H}\UA}$ is also a strictly continuous
function of $\y$. Using again Theorem \ref{af1} we see that $\bhat\tens b$ is
affiliated with $\K{H}\tens\UA$. It is known (\cite[Proposition A.1]{axb})
that if a tensor product of normal operators is affiliated with a tensor
product of $\cst$-algebras then the factor operators
are affiliated with the factor $\cst$-algebras. Thus 
$\bhat\tens b\aff(\K{H}\tens\UA)$ implies that $b\aff\UA$.

Similarly using the strictly continuous family
\[
\Yq\ni\y\longmapsto V(\y)=\chi(\y,a)\in\M{\K{\CC}\tens\K{H}}
\]
we arrive at the conclusion that $\bigl(\chi(\y,a)\bigr)_{\y\in\Yq}$ is a
strictly continuous family of unitary elements of $\M{\UA}$. Consequently $a$
and $a^{-1}$ are affiliated with $\UA$. 

Since, at the same time, $a$, $a^{-1}$ and $b$ generate $\UB$, we see that the
identity mapping on $\UB$ is a morphism from $\UB$ o $\UA$. In other words
$\UB\subset\M{\UA}$ and $\UB\UA$ is a dense subset of $\UA$. This concludes
the proof of the equality $\UA=\UB$.
\qed

\subsection{Quantum group structure}

Having constructed a multiplicative unitary operator and established its
modularity we can proceed with construction of new quantum ``$az+b$'' groups. 

The algebra $\UA$ carries a comultiplication $\del\in\Mor{\UA}{\UA\tens\UA}$.
This is a coassociative morphism given by
\[
\del(c)=W(c\tens I)W^*.
\]
The next ingredient we are going to examine is the scaling group. It is the
one parameter group $(\tau_t)_{t\in\RR}$ of automorphisms of $\UA$ given by
\[
\tau_t(c)=Q^{2it}cQ^{-2it}.
\]
It is easy to check that
\[
\begin{array}{r@{\;=\;}l@{\smallskip}}
\tau_t(a)&a,\\
\tau_t(b)&q^{2it}b
\end{array}
\]
for all $t\in\RR$. Then the well know formula $\Wtil^*=W^{\top\tens R}$ helps
in determining the unitary antipode $R$. It is the $*$-antiautomorphism of
$c\mapsto c^R$ of $\UA$ given on generators as
\[
\begin{array}{r@{\;=\;}l@{\smallskip}}
a^R&a^{-1},\\
b^R&-qa^{-1}b.
\end{array}
\]
The formula for the polar decomposition of the antipode gives now
\[
\begin{array}{r@{\;=\;}l@{\smallskip}}
\kappa(a)&a^{-1},\\
\kappa(b)&-a^{-1}b.
\end{array}
\]
All these formulae agree with those derived in the Hopf $*$-algebra framework
in Subsection \ref{alg_azb} with $\lambda=q^2$. This is why we call our 
quantum group a quantum
``$az+b$'' group for the deformation parameter $q$.

The multiplicative unitary also provides information about the reduced dual 
$(\widehat{\UA},\widehat{\del})$ of
our quantum group. In the case of our quantum ``$az+b$'' groups the situation
is described by the following proposition (recall that $\ahat=b^{-1}$ and
$\bhat=b^{-1}a$).

\begin{Prop}
The operators $\ahat$ and $\bhat$ are affiliated with $\widehat{\UA}$. There
exists a $\cst$-isomorphism $\Psi\colon\UA\to\widehat{\UA}$ such that
$\Psi(a)=\ahat$, $\Psi(b)=\bhat$ and 
\[
\widehat{\del}\comp\Psi=\sigma(\Psi\tens\Psi)\comp\del,
\]
where $\sigma$ is the flip on $\widehat{\UA}\tens\widehat{\UA}$.
\end{Prop}

The proof of this proposition is identical to that of \cite[Theorem 7.1]{azb}.
We shall only point out
that we can avoid choosing a special representation of $\UB$ (as in
\cite{azb}) by noticing that 
\[
\bigl\{\chi(\y,\ahat)\chi(\bhat,\y'):\:\y,\y'\in\Yq\bigr\}
=\bigl\{\chi(b,\y)\chi(\y',a):\:\y,\y'\in\Yq\bigr\},
\]
so that the multiplicities of the pairs $(a,b)$ and $(\ahat,\bhat)$ are the
same regardless of the chosen representation 
(cf.~\cite[Twierdzenie 7.28]{phd}).

According to the results of \cite{haar}, for a quantum group arising from a
modular multiplicative unitary $W$ with associated operators $\Wtil$, $Q$ and
$\Qhat$, the weight
\[
h(c)=\Tr{\Qhat c\Qhat}
\]
is right invariant and if it is {\em locally finite}\/ (densely defined) then
it is the right Haar measure. It turns out that in our case this weight is 
locally finite. More precisely one finds that for $c=g(a)f(a)$ we
have
\[
h(c^*c)=\Int_{\Yq}|g(\y)|^2\,d\Mu(\y)\Int_{\Ybar}|f(\y)|^2|\y|^2\,d\Mu(\y).
\]
In particular $(\UA,\del)$ together with $\bigl(\kappa,(\tau_t)_{t\in\RR},R,h)$
is a {\em weighted Hopf\/ $\cst$-algebra}\/ as defined in 
\cite[Definition 1.5]{mnw}. It is well known that weighted Hopf $\cst$-algebras
are the same objects as {\em reduced\/ $\cst$-algebraic quantum groups}\/ 
defined in \cite[Definition 4.1]{kv}. The left Haar measure is explicitly
given as $h^{\text{\tiny L}}=h\comp R$. In other words 
$(\UA,\del,h^{\text{\tiny L}},h)$ is a reduced $\cst$-algebraic quantum group.

There is evidence that our quantum ``$az+b$'' groups are amenable in the sense
that their reduced dual coincides with the universal one. It has been proved
for the case of the deformation parameter assuming values between $0$ and $1$
(cf.~\cite{azb,puso}).

\section*{Acknowledgements}

The author wishes to express his gratitude to professor S.L.~Woronowicz who
inspired and helped in developing this work. He also wants to thank professors
M.~Bo\.zejko and W.~Pusz whose many comments and remarks have been invaluable.
This paper was prepared during author's stay at the Mathematisches Institut of
the Westf\"alische Wilhelms-Universit\"at in M\"unster. He would like to thank
professor J.~Cuntz for warm hospitality and perfect atmosphere for scientific
activity.


\begin{thebibliography}{XX}
\void{
\bibitem{baaj}
{\sc Baaj, S.~\& Skandalis, G.} Unitaries muliplicatifs et dualit\'{e}
pour les poiduits crois\'{e}s de $\cst$-alg\'{e}bres,
{\em Ann.~scient.~\'{E}c.~Norm.~Sup.}\/ \textbf{26} $4^{\text{e}}$ s\'{e}rie 
(1993), 425--488.
     }
\bibitem{ncgnt}
{\sc Kruszy\'{n}ski, P.~\& Woronowicz, S.L.} A non-commutative
Gelfand-Naimark theorem, {\em J.~Op.~Theory}\/ \textbf{8} (1982), 361--389.
\bibitem{kv}
{\sc Kustermans J.~\& Vaes, S.} Locally Compact quantum groups, {\em
Ann.~Sci.~Ec. Norm.~Sup.}\/ \textbf{33} No.~4 (2000), 837--934.
\bibitem{mnw}
{\sc Masuda, T., Nakagami, Y.~\& Woronowicz, S.L.} A $\cst$-algebraic
framework for quantum groups, {\em Int.~J.~Math.}\/ \textbf{14} No.~9 (2003),
903--1001.
\bibitem{puso}
{\sc W.~Pusz \& P.M.~So\l{}tan,} Functional form of unitary
representations of the quantum ``$az+b$'' group, {\em Rep.~Math.Phys.}\/
\textbf{52} No.~2 (2003) 309--319.
\bibitem{gl}
{\sc Pusz, W.~\& Woronowicz, S.L.} A quantum $GL(2,\CC)$ group at roots of
unity, {\em Rep.~Math.~Phys.}\/ \textbf{47}, No.~3 (2001), 431--462.
\bibitem{phd}
{\sc So\l{}tan, P.M.} Nowe deformacje grupy afinicznych przekszta\l{}ce\'n
p\l{}aszczyzny, PhD thesis Warsaw University 2003.
\bibitem{modmu}
{\sc So{\l}tan, P.M.~\& Woronowicz, S.L.} A remark on manageable
multiplicative unitaries, {\em Lett.~Math.~Phys.}\/ \textbf{57} (2001), 
239--252.
\bibitem{VDH}
{\sc Van Daele, A.} The Haar measure on some locally compact
quantum groups, preprint OA/0109004.
\bibitem{dualth}
{\sc Woronowicz, S.L.} Duality in the $\cst$-algebra theory,
{\em Proc.~Int.~Congr.~Math.~Warsaw 1983,}\/ PWN Polish Scientific Publishers,
Warsaw, 1347--1356.
\bibitem{unbo}
{\sc Woronowicz, S.L.} Unbounded elements affiliated with $\cst$-algebras
and non-compact quantum groups, {\em Comm.~Math.~Phys.}\/ \textbf{136} (1991),
399--432.
\bibitem{OpEq}
{\sc Woronowicz, S.L.} Operator equalities related to the quantum $E(2)$
group, {\em Comm. Math.~Phys.}\/ \textbf{144} (1992), 417--428.
\bibitem{e2}
{\sc Woronowicz, S.L.} Quantum $E(2)$ group and its Pontryagin dual,
{\em Lett.~Math. Phys.}\/ \textbf{23} (1991), 251--263.
\bibitem{gen}
{\sc Woronowicz, S.L.} $\cst$-algebras generated by unbounded elements,
{\em Rev.~Math. Phys.}\/ \textbf{7}, No.~3 (1995), 481--521.
\bibitem{mu}
{\sc Woronowicz, S.L.} From multiplicative unitaries to quantum groups,
{\em Int.~J.~Math.}\/ \textbf{7}, No.~1 (1996), 127--149.
\bibitem{qef}
{\sc Woronowicz, S.L.} Quantum exponential function, {\em Rev.~Math.~Phys.}\/ 
\textbf{12}, No.~6, (2000), 873--920.
\bibitem{azb}
{\sc Woronowicz, S.L.} Quantum `$az+b$' group on complex plane, {\em 
Int.~J.~Math.}\/ \textbf{12}, No.~4 (2001), 461--503.
\bibitem{haar}
{\sc Woronowicz, S.L.} Haar weight on some quantum groups, preprint University
of Warsaw (2003).
\bibitem{axb}
{\sc Woronowicz, S.L.~\& Zakrzewski, S.} Quantum `$ax+b$' group, {\em 
Rev.~Math. Phys.}\/ \textbf{14}, Nos.~7 \& 8 (2002), 797--828.
\end{thebibliography}
\end{document}